\documentclass[11pt]{article}

\usepackage[utf8]{inputenc} %
\usepackage[T1]{fontenc}    %
\usepackage[table]{xcolor}  %
\definecolor{forest}  {rgb}{0,.4,0} 
\definecolor{midnight}  {rgb}{0,0,.5} 
\usepackage[pdftex,plainpages=false,pdfpagelabels,colorlinks=true,urlcolor=midnight,citecolor=forest]{hyperref}
\usepackage{graphicx}
\usepackage{url}            %
\usepackage{booktabs}       %
\usepackage{amsfonts}       %
\usepackage{nicefrac}       %
\usepackage{xfrac}
\usepackage{microtype}      %
\usepackage{lipsum}		%
\usepackage{natbib}
\usepackage{doi}

\usepackage{multirow}
\usepackage[utf8]{inputenc} %
\usepackage[T1]{fontenc}    %
\usepackage{hyperref}       %
\usepackage{url}            %
\usepackage{booktabs}       %
\usepackage{amsfonts}       %
\usepackage{nicefrac}       %
\usepackage{microtype}      %
\usepackage{xcolor}         %
\usepackage[labelfont=bf]{caption}
\usepackage{subcaption}
\usepackage{mathtools}
\usepackage{svg}
\usepackage{listings}
\usepackage[linesnumbered,boxed,algoruled,noline,noend]{algorithm2e}
\usepackage{caption} %
\usepackage{pdflscape}
\let\oldnl\nl%
\newcommand{\nonl}{\renewcommand{\nl}{\let\nl\oldnl}}%

\SetCommentSty{mycommfont}

\DeclareMathOperator* {\argmin} {argmin~}

\newcommand{\revone}[1]{#1}
\newcommand{\revtwo}[1]{#1}

\usepackage[capitalize,noabbrev]{cleveref}
\usepackage[letterpaper,top=2cm,bottom=2cm,left=3cm,right=3cm,marginparwidth=1.75cm]{geometry}
\usepackage{bm}

\usepackage{tikz}
\usepackage{pgfplots}
\usetikzlibrary{positioning,fit,backgrounds}
\definecolor{bblue}{HTML}{4F81BD}
\definecolor{rred}{HTML}{C0504D}
\definecolor{ggreen}{HTML}{9BBB59}
\definecolor{ppurple}{HTML}{9F4C7C}
\definecolor{yyellow}{HTML}{B8BD4F}

\newcommand{\R}{\mathbb{R}}
\newcommand{\C}{\mathbb{C}}
\newcommand{\bx}{\boldsymbol{x}}
\newcommand{\bA}{\boldsymbol{A}}
\newcommand{\cF}{\mathcal{F}}

\usepackage{chngcntr}
\usepackage{apptools}
\usepackage{lineno} %

\newtheorem{theorem}{Theorem}[section]
\AtAppendix{\counterwithin{theorem}{section}}

\title{
Hardware Acceleration for HPS Algorithms \\ in Two and Three Dimensions}

\usepackage{authblk}

\author[1]{Owen Melia}
\author[1,2]{Daniel Fortunato}
\author[3]{Jeremy Hoskins}
\author[3,4,5]{Rebecca Willett}
\affil[1]{Center for Computational Mathematics, Flatiron Institute, USA}
\affil[2]{Center for Computational Biology, Flatiron Institute, USA}
\affil[3]{Computational and Applied Mathematics, Department of Statistics, University of Chicago, USA}
\affil[4]{Department of Computer Science, University of Chicago, USA}
\affil[5]{Data Science Institute, University of Chicago, USA}
\date{}                   
\begin{document}
\maketitle

\begin{abstract}
We provide a flexible, open-source framework for hardware acceleration, namely massively-parallel execution on general-purpose graphics processing units (GPUs), applied to the hierarchical Poincar\'e--Steklov (HPS) family of algorithms for building fast direct solvers for linear elliptic partial differential equations.
To take full advantage of the power of hardware acceleration, we propose two variants of HPS algorithms to improve performance on two- and three-dimensional problems. 
In the two-dimensional setting, we introduce a novel recomputation strategy that minimizes costly data transfers to and from the GPU; 
in three dimensions, we modify and extend the adaptive discretization technique of \citet{geldermans_adaptive_2019} to greatly reduce peak memory usage.
We provide an open-source implementation of these methods written in JAX, a high-level accelerated linear algebra package, which allows for the first integration of a high-order fast direct solver with automatic differentiation tools. 
We conclude with extensive numerical examples showing our methods are fast and accurate on two- and three-dimensional problems.

\end{abstract}

\section{Introduction}
\label{sec:intro}

Many problems in scientific computing require solving systems of linear, elliptic partial differential equations (PDEs).
Such PDEs can accurately model a variety of physics, such as wave propagation, electrostatics, and diffusion phenomena. 
Because analytical solutions of these equations are often unknown, the task of computing numerical solutions has been an area of active research for hundreds of years. 
Today, there are a myriad of numerical solution methods available, and many are tailored to particular classes of equations or to specific use cases. 
We are most interested in designing methods for settings such as inverse or control problems, where the PDE implicitly defines some functional which, along with its gradient, is evaluated sequentially hundreds or thousands of times in the inner loop of an iterative algorithm. 

In these settings, fast direct solvers \citep{martinsson_fast_2019} are a compelling choice. 
These solvers are able to rapidly compute a high-accuracy solution operator and can rapidly evaluate the solution given new data by applying the solution operator, often at the cost of a few matrix-vector multiplications.
Fast direct solvers are also preferable for certain PDEs with oscillatory solutions \citep{gillman_spectrally_2015}, especially ones modeling wave propagation, as they do not incur a data-dependent iteration complexity cost associated with iterative solvers, which can be quite large \citep{ernst_why_2012}.

Recently, scientific computing has undergone a paradigm shift with the advent of general-purpose hardware accelerators, such as GPUs. 
These hardware accelerators allow for massively parallel computation---they have thousands of processor cores on a single chip---but have strict memory constraints, a resource profile very different from standard multicore CPU architectures. 
This paper explores hardware acceleration of fast direct solvers and introduces new methods to facilitate this acceleration.

In particular, we focus on the hierarchical Poincar\'e--Steklov family of algorithms \citep{martinsson_direct_2013,gillman_direct_2014,gillman_spectrally_2015}, a class of direct solution methods for variable-coefficient elliptic PDEs. 
These methods are characterized by a nested dissection approach combined with a high-order composite spectral discretization. 
We identify the algorithmic structure of these algorithms which makes them amenable to GPU acceleration and introduce new techniques for reducing the memory footprint of these methods. For two-dimensional problems, we introduce a novel recomputation strategy that minimizes data transfer between the GPU and host memory. In three dimensions, we use an adaptive discretization method, which greatly reduces the algorithm's peak memory complexity.
Our numerical examples show these ideas are useful in challenging applied settings such as wave propagation, inverse problems, and molecular biology simulations.

We focus on solving linear, elliptic partial differential equations of the form
\begin{align}
    \label{problem:elliptic_BVP}
    \mathcal{L}u(x) = f(x), & \qquad x \in \Omega,  \\
    \label{problem:dirichlet_data}
        u(x) = g(x), & \qquad x \in \partial \Omega.
\end{align}
In \cref{problem:elliptic_BVP}, $\mathcal{L}$
is a linear, elliptic, second-order partial differential operator with spatially-varying coefficient functions, and $\Omega$ is a square $\subset \R^2$ or cube $\subset \R^3$. \cref{problem:dirichlet_data} specifies Dirichlet boundary data, but our methods can also solve problems with Robin or Neumann boundary data. 
In these problems, we assume we can evaluate the differential operator $\mathcal{L}$, the source $f$, and the boundary data $g$ at a set of discretization points of our choosing. We represent the solution $u$ by its restriction to the same set of discretization points and rely on high-order polynomial interpolation to evaluate $u$ away from the discretization points.
In this paper, we refer to vectors with bold lowercase symbols such as $\boldsymbol{f}$ and matrices with  bold uppercase symbols such as $\bA$.
We use $x$ for the spatial variable, and when we want to indicate Cartesian coordinates, we use 
$(x_1,x_2) \in \R^2$  and $(x_1,x_2,x_3) \in \R^3$. We use the subscript $u_n$ to denote the outward-pointing boundary normal derivative of a function, and we use $\Delta$ to denote the Laplace operator, the sum of second derivatives in each dimension.

\subsection{Paper outline and contributions}
In \cref{sec:related_work}, we discuss related work, including algorithmic development for fast direct solvers and GPU-specific optimizations. 
In \cref{sec:hps_algos}, we give an overview of the hierarchical Poincar\'e--Steklov method and discuss the potential for massively parallel implementations of the algorithm. 
In the rest of the paper, we make the following contributions:
\begin{itemize}
\item We optimize data transfer patterns to accelerate our method applied to two-dimensional problems (\cref{sec:hardware_acceleration_2D}).
\item To alleviate peak memory requirements in three-dimensional problems, we extend the two-dimensional adaptive method of \cite{geldermans_adaptive_2019} to three dimensions, develop the first adaptive 3D GPU-compatible HPS implementation, and provide numerical examples to demonstrate memory and accuracy tradeoffs (\cref{sec:hardware_acceleration_3D}).
\item We provide a range of numerical examples illustrating the application of our method, focusing on two settings: high-wavenumber scattering problems and the linearized Poisson--Boltzmann equation (\cref{sec:2D_examples,sec:3D_examples}).
\item We show our proposed algorithm and implementation can be combined easily with standard automatic differentiation software, which makes it particularly amenable to application in optimization, inverse problems, and machine learning contexts (\cref{sec:2D_examples}).
\item We make our JAX-based implementation publicly available at \url{https://github.com/meliao/jaxhps}.
\end{itemize}

\section{Related work}
\label{sec:related_work}

HPS algorithms are built on two conceptual building blocks: composite high-order spectral collocation methods
\citep{kopriva_staggered-grid_1998,pfeiffer_multidomain_2003,yang_multidomain_2000}, and nested dissection of the computational domain \citep{george_nested_1973}.
Composite spectral collocation methods are those which separate the computational domain into a set of disjoint elements, and use a high-order spectral collocation scheme to represent the problem and solution separately on each element.
Nested dissection methods break the original problem into a series of subproblems defined on a hierarchy of subdomains. 
The careful ordering of subproblems reduces the overall computational complexity by leveraging knowledge about properties of the solution, i.e. continuity of the solution and its derivative. 
Composite spectral collocation and hierarchical matrix decomposition ideas were combined in integral equation methods \citep{ho_fast_2012} for constant-coefficient PDEs. 

\citet{martinsson_direct_2013} first proposed combining these elements in a fast direct solver for variable-coefficient linear elliptic PDEs. The proposed scheme discretizes and merges Dirichlet-to-Neumann (DtN) operators. 
\citet{gillman_direct_2014} proposed a compression scheme that leverages the structure of these DtN operators to build a solver with $O(n)$ computational complexity for $n$ elements.
To alleviate the instabilities observed when merging DtN operators for Helmholtz problems, \citet{gillman_spectrally_2015} proposed a scheme that merges impedance-to-impedance (ItI) operators instead.
Further analysis for this scheme was provided in \citet{beck_quantitative_2022}.
Modifications for three dimensions have been proposed, including \citet{lucero_lorca_iterative_2024,hao_direct_2016,Kump_Yesypenko_Martinsson_2025}, which all build solvers for three-dimensional Helmholtz problems. To alleviate memory and computational complexity, \citet{lucero_lorca_iterative_2024} use an iterative method at the highest-level subproblems.
In concurrent work to our own, \cite{Kump_Yesypenko_Martinsson_2025} approach three-dimensional problems with uniform discretizations using a hybrid GPU-CPU approach by combining the composite spectral collocation method with a two-level sparse direct solver. 
\citet{Fortunato_Hale_Townsend_2021} use the ultraspherical spectral method to discretize triangular or quadrilateral mesh elements and compute solutions over polygonal domains by merging DtN operators.
\citet{fortunato_high-order_2024} develops a variant of the HPS method which merges DtN and ItI operators to solve PDEs on unstructured meshes of smooth two-dimensional surfaces. 
\revtwo{ \mbox{\citet{beams_parallel_2019}}  develops an implementation of the HPS algorithm targeting parallel shared-memory computer architectures. }

\sloppypar{
There has been significant interest in the GPU acceleration of (low-order) iterative PDE solvers \citep{georgescu_gpu_2013}. 
\revone{Other general-purpose packages, such as MFEM and libCEED, implement high-order iterative solvers with GPU acceleration \mbox{\citep{Abdelfattah_GPU_2021,Kolev_efficient_2021}}.
}
Accelerating these algorithms often requires the rapid application of an extremely sparse system matrix. 
Applying GPU acceleration to direct solvers \citep{abdelfattah_addressing_2022,ghysels_high_2022,li_superlu_dist_2003} requires different techniques; 
the literature has mostly focused on sparse direct solvers which do not employ a nested dissection method. 
}

These sparse direct solvers often have much higher peak memory requirements and heterogeneous computation profiles when compared with iterative PDE solvers.

Other solvers have been designed directly for GPU acceleration. 
\citet{yesypenko_gpu_2024,yesypenko_slablu_2024} designed a composite high-order spectral collocation method and associated sparse direct solver with highly heterogeneous computation patterns, which eases GPU acceleration. 
This method can solve variable-coefficient 2D problems very quickly; 
our method solves similar problems, but can also handle three-dimensional problems and interface with automatic differentiation.
Developing a GPU-compatible implementation of the solver in \citet{yesypenko_slablu_2024}, as well as the implementation of our method, has been greatly eased by the advent of high-performance hardware-accelerated linear algebra frameworks popularized by deep learning, such as PyTorch \citep{ansel_pytorch_2024} and JAX \citep{jax2018github}. These frameworks are highly efficient for batched linear algebra tasks and implement automatic differentiation capabilities.
Both are high-level packages that sit on top of the XLA compiler \citep{Leary_Wang_2017}, which compiles and launches optimized kernels that execute on general-purpose GPUs.
There has also been work creating automatic differentiation compatible PDE solvers in JAX, 
tailored for problems such as synchrotron simulation \citep{diao_synax_nodate}, 
computational mechanics \citep{xue_jax-fem_2023}, 
and ordinary differential equations \citep{kidger_neural_2021}.

\section{Introduction to HPS methods}
\label{sec:hps_algos}
\newcommand{\nleaves}{n_{\text{leaves}}}
In this section, we provide an overview of the HPS algorithms used in the paper, with a particular eye on their computational structure and the possibilities for GPU acceleration.\footnote{
    \revtwo{For an instructional introduction to these methods, we refer the interested reader to \mbox{\citet{martinsson_hierarchical_2015,gillman_spectrally_2015,martinsson_fast_2019}}.}
    }
Full algorithms are available in \cref{appendix:2D_DtN,appendix:2D_ItI,appendix:3D_DtN_uniform}.
We use different variants of this algorithm for merging different types of Poincar\'e--Steklov operators. 
A Poincar\'e--Steklov operator $T: g \mapsto h$ maps from one type of boundary data to another. Take, for example, a Dirichlet-to-Neumann (DtN) operator, which maps from Dirichlet data $g$ on the boundary of $\Omega$ to Neumann data $h$ on the same boundary:
\begin{linenomath}
\begin{linenomath}\begin{align*}
    g &= u|_{\partial \Omega}, \\
    h &= u_n|_{\partial \Omega},
\end{align*}\end{linenomath}
\end{linenomath}
where $u$ satisfies \cref{problem:elliptic_BVP}.
Another example commonly used is an impedance-to-impedance (ItI) operator, which maps ``incoming'' \revtwo{impedance data $g$} to ``outgoing'' impedance data \revtwo{$h$} \citep{gillman_spectrally_2015}:
\begin{linenomath}\begin{align*}
    g &= u_n + i\eta u|_{\partial \Omega}, \\
    h &= u_n - i\eta u|_{\partial \Omega},
\end{align*}\end{linenomath}
where $u$ satisfies \cref{problem:elliptic_BVP}.
These Poincar\'e--Steklov operators are linear operators, and we work with their discretization $\boldsymbol{T}$. Throughout the algorithm, we also work with $\boldsymbol{g}$ and $\boldsymbol{h}$, vectors of incoming and outgoing boundary data evaluated at a set of discretization points.

It is important to remember that because we are solving a linear partial differential equation, we can decompose the solution $u(x) $ into a particular solution $v(x)$ and homogeneous solution $w(x)$ where $u(x) = v(x) + w(x)$. 
The particular solution $v(x)$ satisfies
\begin{linenomath}\begin{align*}
    \left\{
    \begin{aligned}
        \mathcal{L} v(x) &= f(x), &&\quad x \in \Omega, \\
        v(x) &= 0, &&\quad x \in \partial \Omega,
    \end{aligned}
        \right.
\end{align*}\end{linenomath}
and the homogeneous solution $w(x)$ satisfies
\begin{linenomath}\begin{align*}
    \left\{
    \begin{aligned}
        \mathcal{L} w(x) &= 0, &&\quad x \in \Omega, \\
        w(x) &= g(x), &&\quad x \in \partial \Omega.
    \end{aligned}
        \right.
\end{align*}\end{linenomath}

\subsection{Discretization via composite high-order spectral collocation}
\label{sec:discretization}
To numerically solve \cref{problem:elliptic_BVP}, \revtwo{a discretization is needed to represent} $\mathcal{L}$, $f$, and $g$ in some finite-dimensional basis. \revtwo{HPS methods perform this discretization} in two steps: a recursive partition of the domain $\Omega$ and a high-order spectral collocation scheme. 
The first step is to recursively partition $\Omega$ using a quadtree or octree structure down to a user-specified maximum depth $L$. 

We use $\nleaves$ to denote the number of patches at the finest level of the spatial partition.
In this section, we consider uniform discretization trees, so each tree with depth $L$ will have $\nleaves=2^{d L}$ elements at the lowest level in dimension $d = 2,3$. In \cref{sec:hardware_acceleration_3D}, we consider more general discretization trees.
At times, it will be useful to describe the progress of the algorithm using language to describe the trees representing the spatial partition. 
To that end, we will sometimes refer to the elements as \emph{nodes} and elements at the lowest level of the tree as \emph{leaves}. The element at the highest level of the tree, which represents the entire computational domain, is sometimes called the \emph{root}.

\begin{figure}
    \centering
    \begin{subfigure}[b]{0.3\textwidth}
        \centering
        \includeinkscape[width=\linewidth]{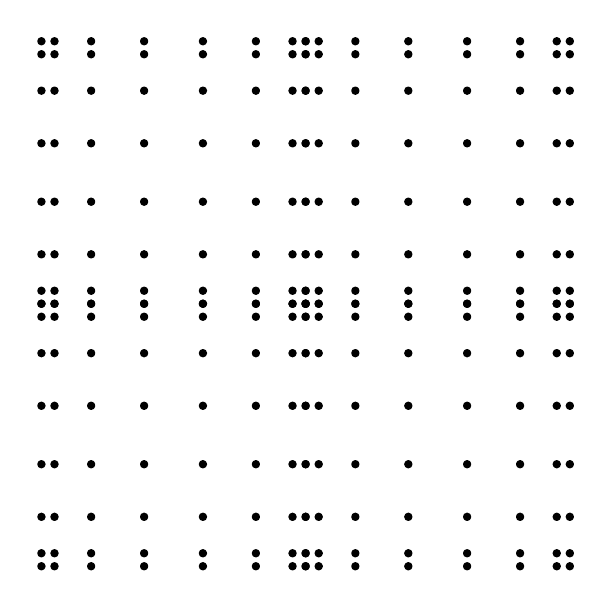}
        \caption{}
        \label{fig:chebyshev_points}
    \end{subfigure}%
    \begin{subfigure}[b]{0.3\textwidth}
        \centering
        \includeinkscape[width=\linewidth]{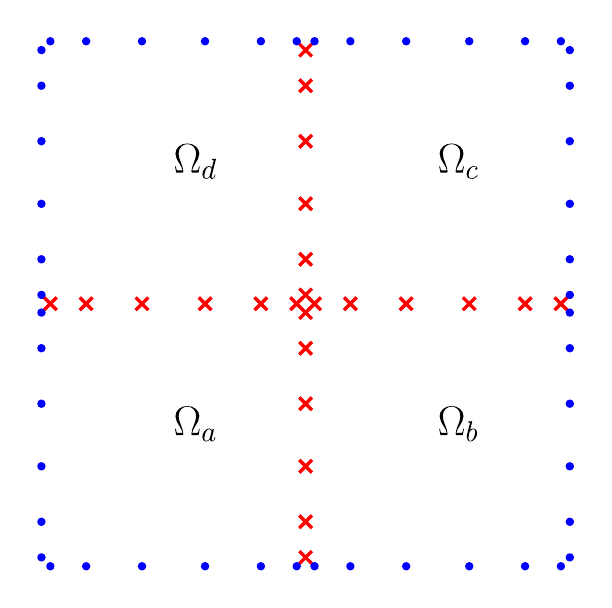}
        \caption{}
        \label{fig:quad_merge_int_ext_points}
    \end{subfigure}
    \caption{Visualizing the high-order composite spectral collocation scheme for a simple two-dimensional problem. \cref{fig:chebyshev_points} shows the Chebyshev points on a two-dimensional problem with polynomial order $p=8$ and $L=1$ level of refinement. \cref{fig:quad_merge_int_ext_points} shows the Gauss--Lobatto points discretizing the boundaries of the leaves using order $q=6$. When merging nodes together, it is important to distinguish between the \emph{exterior} boundary points, drawn with blue dots, and the \emph{interior} boundary points, drawn with red $\times$'s. 
    }
    \label{fig:disc_points}
\end{figure}

\revtwo{Each leaf is discretized} using a tensor product of Chebyshev--Lobatto points, with user-specified order $p$. 
This requires $p^d$ points per leaf in dimension $d=2,3$.
\revtwo{HPS methods typically call for an order-$q$ Gauss--Legendre quadrature rule to} represent the boundary of each leaf. 
For simplicity and stability, we always use $q = p-2$\revtwo{, following \mbox{\citet{gillman_spectrally_2015}}.}
In two and three dimensions, there are $(2 d)q^{d-1}$ boundary discretization points per leaf. We show the interior and boundary points in \cref{fig:disc_points}. The resulting discretization has $N= \nleaves p^d = (2^Lp)^d$ interior discretization points.

\subsection{Local solve stage}
\label{sec:local_solve}
\revtwo{The first step is to} discretize the differential operator restricted to that leaf on the $p^d$ Chebyshev discretization points. We call the resulting matrix $\boldsymbol{L}^{(i)}$ for leaf $i$. This matrix is a combination of Chebyshev spectral differentiation matrices \citep{trefethen_spectral_2000} and evaluations of the spatially varying coefficient functions. 
\revtwo{The} source function $f$ \revtwo{is discretized} on the same points, \revtwo{and we} call the resulting vector $\boldsymbol{f}^{(i)}$.
At this point, \revtwo{the HPS algorithm solves} a local boundary-value problem on each patch. 
\revtwo{Standard practice \mbox{\citep{gillman_spectrally_2015,fortunato_high-order_2024}} is to} solve this problem using a ``boundary bordering'' technique which enforces the differential operator on the Chebyshev nodes interior to the leaf, and enforces a Dirichlet or impedance 
boundary condition on the Chebyshev nodes on the leaf's boundary. 
At this point in the algorithm, the correct boundary values to enforce at each leaf \revtwo{are unknown,} so $\boldsymbol{L}^{(i)}$ and $\boldsymbol{f}^{(i)}$ are used to precompute a solution map for the local problem.
To precompute this local solution map,\revtwo{ the algorithm constructs} a matrix $\boldsymbol{Y}^{(i)}$ which maps from any boundary data $\boldsymbol{g}^{(i)}$ to the corresponding homogeneous solution on the Chebyshev nodes. 
\revtwo{The algorithm also computes} $\boldsymbol{v}^{(i)}$, a vector evaluating the particular solution on the Chebyshev nodes. 
\revtwo{Both $\boldsymbol{Y}^{(i)}$ and $\boldsymbol{v}^{(i)}$ can be expressed simply in terms of the input data; lines \mbox{\ref{eqn:def_Y_DtN}} and \mbox{\ref{eqn:dev_v_DtN}} in \mbox{\cref{alg:2D_DtN_local}} in the Appendix give these expressions.}
For each leaf, \revtwo{the algorithm also constructs} a Poincar\'e--Steklov matrix $\boldsymbol{T}^{(i)}$ and a vector of outgoing data $\boldsymbol{h}^{(i)}$.
Computing $\boldsymbol{T}^{(i)}$ and $\boldsymbol{h}^{(i)}$ only requires multiplying $\boldsymbol{Y}^{(i)}$ and $\boldsymbol{v}^{(i)}$ with a fixed, precomputed operator which composes interpolation matrices from Chebyshev to Gauss--Legendre discretization points and spectral Chebyshev differentiation matrices.

\cref{alg:local_solve} shows that this stage of the algorithm is a long loop over linear algebra operations. The size of these operations is controlled by the polynomial order $p$, and we find these operations are efficient for the orders $p \leq 16$ considered in this work.
The units of work inside the loop are \emph{embarrassingly parallel}, meaning that one iteration does not depend on the output of any other iterations. 
Furthermore, because we hold $p$ constant on all leaves, all of the linear algebra operations are homogeneous, meaning they all operate on the same sizes of matrices.
This computational structure facilitates GPU acceleration by batching and parallelizing local solves.
We give full details describing the local solve stage in \cref{alg:2D_DtN_local,alg:2D_ItI_local} in the Appendix. 

\DontPrintSemicolon
\SetArgSty{textnormal}
\begin{algorithm}[h!]
    \caption{
       Local solve stage. Full details are available in \cref{alg:2D_DtN_local,alg:2D_ItI_local}.
    }
    \label{alg:local_solve}
    \KwIn{Discretized differential operators $\{\boldsymbol{L}^{(i)} \}_{i=1}^{\nleaves}$; 
    discretized source functions $\{ \boldsymbol{f}^{(i)} \}_{i=1}^{\nleaves}$; 
    precomputed interpolation and differentiation matrices}
    \For(){Leaf $i=1, \hdots , \nleaves$} {
        Perform boundary bordering to $\boldsymbol{L}^{(i)}$ \\
        Invert the resulting matrix \tcp*{Main computational work}
        Construct $\boldsymbol{Y}^{(i)}$, the interior solution matrix \\
        Construct $\boldsymbol{T}^{(i)}$, the Poincar\'e--Steklov matrix \\
        Construct $\boldsymbol{v}^{(i)}$, the leaf-level particular solution \\
        Construct $\boldsymbol{h}^{(i)}$, the outgoing boundary data \\
    }
    \KwResult{
    Poincar\'e--Steklov matrices $\{ \boldsymbol{T}^{(i)} \}_{i=1}^{\nleaves}$; 
    outgoing boundary data $\{ \boldsymbol{h}^{(i)} \}_{i=1}^{\nleaves}$; 
    interior solution matrices $\{ \boldsymbol{Y}^{(i)} \}_{i=1}^{\nleaves}$;
    leaf-level particular solutions $\{ \boldsymbol{v}^{(i)} \}_{i=1}^{\nleaves}$ }
\end{algorithm}

\subsection{Merge stage}
\label{sec:hps_algos_merge}
After computing $\boldsymbol{T}^{(i)}$ and $\boldsymbol{h}^{(i)}$ for each leaf in the local solve stage, \revtwo{the HPS algorithm begins} merging nodes of the tree together. 
This process creates a hierarchy of solution operators which will later be used to propagate the boundary data on $\partial \Omega$ to the boundary of each leaf.
For 2D problems, \revtwo{our implementation merges} nodes four at a time, and for 3D problems, \revtwo{our implementation merges} nodes eight at a time. 
Suppose \revtwo{the algorithm is} merging a set of nodes $\{a, b, \hdots \}$ which all share parent node $j$. \revtwo{The algorithm has} access to the following data:
\begin{itemize}
    \item $\{ \boldsymbol{T}^{(a)}, \boldsymbol{T}^{(b)}, \hdots \}$, the Poincar\'e--Steklov matrices of the nodes being merged.
    \item $\{ \boldsymbol{h}^{(a)}, \boldsymbol{h}^{(b)}, \hdots \}$, the outgoing boundary data due to the particular solution of the nodes being merged.
\end{itemize}

\newcommand{\ext}{\textrm{ext}}
\newcommand{\intt}{\textrm{int}}
\newcommand{\child}{\textrm{child}}
At this point, it is helpful to distinguish between vectors that are defined along the \emph{exterior} of the patches being merged and vectors that are defined along the \emph{interior} of the patches being merged. We indicate these vectors with subscripts $ext$ and $int$, respectively. 
See \cref{fig:quad_merge_int_ext_points} for a diagram of the interior and exterior points in a 2D merge operation. 
\revtwo{The goal of the merge operation is} to precompute a solution operator which propagates the information from the exterior boundary points to the interior boundary points. 
This solution operator takes the form $\boldsymbol{g}^{(j)}_{\ext} \mapsto \boldsymbol{S}^{(j)} \boldsymbol{g}^{(j)}_{\ext} + \tilde{\boldsymbol{g}}^{(j)}$. 
In this equation, $ \boldsymbol{S}^{(j)} \boldsymbol{g}^{(j)}_{\ext}$ evaluates the homogeneous solution on the interior boundary points, and $\tilde{\boldsymbol{g}}^{(j)}$ evaluates the particular solution on the interior boundary points. 
As in \cref{sec:local_solve}, the boundary data $\boldsymbol{g}^{(j)}_{\ext} $ is not available at this stage of the algorithm, but \revtwo{it is possible to} precompute the other parts of the solution operator. To that end, each merge operation \revtwo{will compute}:
\begin{itemize}
    \item $\boldsymbol{S}^{(j)}$, the propagation operator for node $j$, which maps incoming homogeneous boundary data from the exterior boundary points to the interior boundary points.
    \item $\tilde{\boldsymbol{g}}^{(j)}$, the incoming boundary data due to the particular solution evaluated at the interior boundary points.
    \item $\boldsymbol{T}^{(j)}$, the Poincar\'e--Steklov matrix for node $j$.
    \item $\boldsymbol{h}^{(j)}$, the outgoing boundary data for node $j$.
\end{itemize}
\revtwo{$\boldsymbol{S}^{(j)}$ and $\tilde{\boldsymbol{g}}^{(j)}$} will be used in the final stage of the HPS method when applying the precomputed solution operator, and \revtwo{$\boldsymbol{T}^{(j)}$ and $\boldsymbol{h}^{(j)}$} will be used in a future merge operation.

To compute the merge outputs, \revtwo{the HPS algorithm sets} up a system of equations for a given $\boldsymbol{g}_{\ext}^{(j)}$ and unknown $\boldsymbol{g}_{\intt}^{(j)}$ and solve using a Schur complement approach.\footnote{
    To the best of our knowledge, the block system \revtwo{associated with the merge four box procedure} for \revtwo{the} 2D DtN \revtwo{method was first documented in \cite{chipman_fast_2024}.} 
    We generalize this presentation to encompass 3D problems and merging ItI matrices as well.
    }
The constraints in this system come from our knowledge that the solution and its derivative will be continuous across merge interfaces.
The resulting system of constraints is a linear system and can be written in a blockwise fashion: 
\begin{align}
    \begin{bmatrix} \boldsymbol{A}\vphantom{\boldsymbol{g}_{\ext}^{(j)}} & \boldsymbol{B} \\
        \boldsymbol{C} & \boldsymbol{D}
    \end{bmatrix} 
    \begin{bmatrix}
        \boldsymbol{g}_{\ext}^{(j)} \\ \boldsymbol{g}_{\intt}^{(j)}
    \end{bmatrix} = \begin{bmatrix}
        \boldsymbol{u}^{(j)}_{\ext} -\boldsymbol{h}^{(\child)}_{\ext} \\
        - \boldsymbol{h}_{\intt}^{(\child)}
    \end{bmatrix}.
    \label{eq:merge_lin_system}
\end{align}

\revtwo{The next step is to} build the matrices $\boldsymbol{A},\boldsymbol{B}, \boldsymbol{C}, \boldsymbol{D}$ blockwise from the child Poincar\'e--Steklov matrices $\{ \boldsymbol{T}^{(a)}, \boldsymbol{T}^{(b)}, \hdots \}$, build $\boldsymbol{h}^{(\child)}_{\ext}$ and $ \boldsymbol{h}_{\intt}^{(\child)}$ from the child outgoing boundary data $\{ \boldsymbol{h}^{(a)}, \boldsymbol{h}^{(b)}, \hdots \}$, and use $\boldsymbol{u}^{(j)}_{\ext}$ to represent the unknown solution evaluated on the exterior boundary points.
\revtwo{For 2D DtN merges, these objects are defined in \mbox{\cref{eq:2D_DtN_h_child_ext,eq:2D_DtN_A,eq:2D_DtN_B,eq:2D_DtN_merge_2,eq:2D_DtN_h_child_int,eq:2D_DtN_C,eq:2D_DtN_D}}. Derivations for the 2D ItI and 3D DtN cases can be found in \mbox{\cref{appendix:2D_ItI,appendix:3D_DtN_uniform}}.}
At this point in the algorithm, $\boldsymbol{u}^{(j)}_{\ext}$ \revtwo{is unknown,} which means \revtwo{one }can not directly invert the linear system to solve for $\boldsymbol{g}_{\ext}^{(j)}$ or $\boldsymbol{g}_{\intt}^{(j)}$. However, \revtwo{one} can use a Schur complement approach to partially solve the system:
\begin{linenomath}\begin{align*}
    \begin{bmatrix} \boldsymbol{A} \vphantom{\boldsymbol{g}_{\ext}^{(j)}} - \boldsymbol{B}\boldsymbol{D}^{-1}\boldsymbol{C} & \boldsymbol{0} \\
        \boldsymbol{D}^{-1}\boldsymbol{C} & \boldsymbol{I}
    \end{bmatrix} 
    \begin{bmatrix}
        \boldsymbol{g}_{\ext}^{(j)} \\ \boldsymbol{g}_{\intt}^{(j)}
    \end{bmatrix} = \begin{bmatrix}
        \boldsymbol{u}^{(j)}_{\ext} -\boldsymbol{h}^{(\child)}_{\ext} - \boldsymbol{B} \boldsymbol{D}^{-1} \boldsymbol{h}_{\intt}^{(\child)} \\
        - \boldsymbol{D}^{-1} \boldsymbol{h}_{\intt}^{(\child)}
    \end{bmatrix}.
\end{align*}\end{linenomath}
Interpreting the rows of this linear system gives us the desired outputs:
\begin{align}
    \boldsymbol{u}^{(j)}_{\ext} &= 
    \underbrace{\left( \boldsymbol{A} - \boldsymbol{B}\boldsymbol{D}^{-1}\boldsymbol{C}  \right) }_{ \boldsymbol{T}^{(j)}} \boldsymbol{g}_{\ext}^{(j)}
    + \underbrace{ \boldsymbol{h}^{(\child)}_{\ext}  - \boldsymbol{B} \boldsymbol{D}^{-1}\boldsymbol{h}_{\intt}^{(\child)}}_{\boldsymbol{h}^{(j)}},
    \label{eq:merge_outputs_1} \\
    \boldsymbol{g}_{\intt}^{(j)} &= \underbrace{-\boldsymbol{D}^{-1}\boldsymbol{C}}_{\boldsymbol{S}^{(j)}} \boldsymbol{g}_{\ext}^{(j)} + \underbrace{ -\boldsymbol{D}^{-1}\boldsymbol{h}_{\intt}^{(\child)}.}_{\tilde{\boldsymbol{g}}^{(j)}} \label{eq:merge_outputs_2}
\end{align}
\cref{alg:merge} gives pseudocode for the merge stage of the HPS algorithm. The majority of the computational work for each merge is inverting $\boldsymbol{D}$, which has size proportional to the number of discretization points along the merge interfaces. 
This matrix is quite small at the lowest levels of the discretization tree and grows as the algorithm proceeds to higher nodes in the tree.
Similar to the local solve stage, each merge operation is dominated by linear algebra work, and the units of work inside the inner loop are embarrassingly parallel. 
This means we can easily use GPU acceleration to parallelize the inner loop of \cref{alg:merge}. The outer loop is iterating over different levels, which means the computation during outer loop iteration $\ell$ depends on the outputs of the previous iteration. Because we only use a moderate number of refinement levels $L < 10$, we find \cref{alg:merge} executes very quickly on the GPU despite this dependency structure. We note that our choice to merge nodes four-to-one and eight-to-one (rather than the standard two-to-one) decreases the length of the outer loop by factors of two and three, respectively.

\newlength{\commentWidth}
\setlength{\commentWidth}{4.5cm}
\newcommand{\atcp}[1]{\tcp*[r]{\makebox[\commentWidth]{#1\hfill}}}

\begin{algorithm}[h!]
    \caption{
       Merge stage. Full details are available in \cref{appendix:2D_DtN_merge,appendix:2D_ItI_merge,appendix:3D_DtN_merge}.
    }
    \label{alg:merge}
    \KwIn{Leaf-level Poincar\'e--Steklov matrices $\{\boldsymbol{T}^{(i)} \}_{i=1}^{\nleaves}$; Leaf-level outgoing boundary data $\{ \boldsymbol{h}^{(i)} \}_{i=1}^{\nleaves}$}
    \For() {Merge level $\ell = L-1,\hdots , 0$} {
        \For() {Node $j$ in level $\ell$}{
            Let $a, b, \hdots$ be the children of node $j$ \\
            Use $\{ \boldsymbol{T}^{(a)} , \boldsymbol{T}^{(b)}, \hdots \}$ to build blocks $\boldsymbol{A}, \boldsymbol{B}, \boldsymbol{C}$, and $\boldsymbol{D}$ \\
            Use $\{ \boldsymbol{h}^{(a)} , \boldsymbol{h}^{(b)}, \hdots \}$ to build $\boldsymbol{h}_{\ext}^{(\child)}$ and $\boldsymbol{h}_{\intt}^{(\child)}$ \\
            Invert $\boldsymbol{D}$  \atcp{ Main computational work}
            Evaluate $\boldsymbol{T}^{(j)}$, $\boldsymbol{h}^{(j)}$, $\boldsymbol{S}^{(j)}$, $\tilde{\boldsymbol{g}}^{(j)}$ \atcp{\cref{eq:merge_outputs_1,eq:merge_outputs_2}}
        }
    }
    \KwResult{Poincar\'e--Steklov matrices $\boldsymbol{T}^{(j)} $ for each node; 
    outgoing boundary data $ \boldsymbol{h}^{(j)} $ for each node; 
    propagation operators $ \boldsymbol{S}^{(j)} $ for each node; 
    incoming particular solution data ${\boldsymbol{\tilde{g}}}^{(j)} $ for each node }
\end{algorithm}

\subsection{Downward pass}
In the final stage of the HPS method, all parts of the structured solution operators mapping $g \mapsto u$ have been computed. 
\revtwo{The HPS algorithm evaluates} this structured solution operator by propagating information down the discretization tree from the boundary of the root to the interior of the leaves. The $\boldsymbol{S}^{(j)}$ matrices propagate the homogeneous boundary data to the merge interfaces, and the $\tilde{\boldsymbol{g}}^{(j)}$ vectors add back in the particular solution \revtwo{(line \mbox{\ref{eq:downpass_propagation}} in \mbox{\cref{alg:down_pass}}.)}
This part of the HPS method is extremely fast, as it only involves matrix-vector products. As with the structure of the merge stage, the iterations of the inner loop are embarrassingly parallel and can be batched on the GPU. We show the pseudocode for this stage in \cref{alg:down_pass}.

\begin{algorithm}[h!]
    \caption{
       Downward pass.
    }
    \label{alg:down_pass}
    \KwIn{Boundary data $\boldsymbol{g}$; 
    propagation operators $ \boldsymbol{S}^{(j)} $ for each node; 
    incoming particular solution data ${\boldsymbol{\tilde{g}}}^{(j)}$ for each node;
    leaf-level interior solution matrices $\{ \boldsymbol{Y}^{(i)} \}_{i=1}^{\nleaves}$; 
    leaf-level particular solutions $\{ \boldsymbol{v}^{(i)} \}_{i=1}^{\nleaves}$}
    \For() {Merge level $\ell = 0, ..., L-1$} {
        \For() {Node $j$ in level $\ell$}{
            Look up $\boldsymbol{S}^{(j)}, \boldsymbol{g}^{(j)}$, and $ \tilde{\boldsymbol{g}}^{(j)}$\\
            $\boldsymbol{g}_{\intt} = \boldsymbol{S}^{(j)} \boldsymbol{g}^{(j)} + \tilde{\boldsymbol{g}}^{(j)}$ \label{eq:downpass_propagation}\\
            Let $a, b, \hdots$ be the children of node $j$ \\
            Concatenate $\boldsymbol{g}_{\intt}$ and $\boldsymbol{g}^{(j)}$ to form $\{ {\boldsymbol{g}}^{(a)} , {\boldsymbol{g}}^{(b)}, \hdots \}$
        }
    }
    \For() {Leaf $i=1,\dots , \nleaves$}{
            $\boldsymbol{u}^{(i)} = \boldsymbol{Y}^{(i)} \boldsymbol{g}^{(i)}+ \boldsymbol{v}^{(i)}$
        }
    \KwResult{Leaf-level solutions on the Chebyshev discretization points $\{\boldsymbol{u}^{(i)} \}_{i=1}^{\nleaves}$ }
\end{algorithm}

\section{Hardware acceleration for HPS methods}
\label{sec:hardware_acceleration}
\sloppypar{
While the algorithms presented in \cref{sec:hps_algos} have attractive computational complexity ($O(p^6\nleaves + p^3\nleaves^{3/2})$ in 2D and $O(p^9\nleaves + p^6\nleaves^{2})$ in 3D) and possibilities for parallel execution, they also incur large memory footprints. 
At each step of the algorithm, dense solution matrices are precomputed and must be stored for future use.
\revone{The outputs of the merge stage dominate the memory complexity of the method. In 2D, the overall memory complexity of storing these matrices is $O\left( p^2 \nleaves L \right)$, and in 3D, the memory complexity is $O\left( p^4 \nleaves 2^L \right)$, with a large prefactor.}
This poses a significant challenge for GPU acceleration as general-purpose GPU architectures have significantly more processor cores per unit of memory than standard multicore CPU nodes. 
While standard multicore compute nodes may have 1TB of available random-access memory (host RAM), high-end GPUs have only 80GB of on-device memory, with slow interconnects between the GPU and host RAM. 
This means that batched linear algebra operations are extremely fast on the GPU, but the overall algorithm is slowed by steps transferring data between the GPU and the host.
Thus, to efficiently accelerate HPS algorithms on the GPU, one must devote significant thought to reducing the memory footprint of these algorithms. 
In this section, we introduce two ideas to reduce this memory footprint in two and three-dimensional problems.
}

\subsection{Recomputation strategies to minimize communication costs}
\label{sec:hardware_acceleration_2D}
CPU-bound implementations of HPS methods in 2D often spend an order of magnitude longer in the local solve stage than in the merge or downward pass stages \citep{fortunato_high-order_2024}. 
This suggests that HPS schemes can be greatly accelerated by placing the local solve stage computation on the GPU alone. Indeed, such savings have been observed in \cite{yesypenko_gpu_2024}. 
We observe that for the lowest levels of the merge stage, each unit of computational work is similarly small, suggesting that the GPU can efficiently accelerate this part of the algorithm, too, provided the algorithm is correctly expressed to leverage its inherent parallelism.
In many GPUs, the interconnect between the host device and the GPU's on-device memory is very slow in relation to the speed at which the thousands of processor cores can process data. 
This means that for large problem sizes, operations transferring precomputed solution operators to and from the GPU are a major impediment to fast execution as they incur a large latency and are often ``blocking'' operations, which require all parallel threads to complete before executing.

For large problem sizes in 2D, it is advantageous to delete some data computed in the early stages of the algorithm and recompute it later when necessary. 
This strategy increases the number of floating point operations on the GPU but minimizes costly data transfers. 
A leaf-level recomputation strategy for implementing HPS algorithms on a GPU is presented in \cref{alg:baseline_recomputation}; a similar method is presented in \cite{yesypenko_gpu_2024}. 
The leaf recomputation strategy avoids transferring the $\{ \boldsymbol{Y}^{(i)} \}_{i=1}^{\nleaves}$ and  $\{ \boldsymbol{v}^{(i)} \}_{i=1}^{\nleaves}$ by performing the local solve stage again at the end of the algorithm. 
Under this recomputation strategy, all of the leaf-level Poincar\'e--Steklov matrices must be transferred to RAM during the local solve stage, and then back to the GPU during the merge stage.

\begin{figure}[h!]
    \begin{subfigure}[b]{0.5\textwidth}
    \centering
    \resizebox{\linewidth}{!}{
    \begin{tikzpicture}[
        level 1/.style={sibling distance=100mm},
        level 2/.style={sibling distance=50mm},
        level 3/.style={sibling distance=25mm},
        every node/.style={circle, draw, minimum size=20pt},
        edge from parent/.style={draw, -},
        box/.style={rectangle, draw=gray, dashed, rounded corners, inner sep=8pt, line width=1.5pt}
    ]
        \node (n1) {1}
            child {
                node (n2) {2}
                child {
                    node (n4) {4}
                    child {
                        node (n8) {8}
                    }
                    child {
                        node (n9) {9}
                    }
                }
                child {
                    node (n5) {5}
                    child {
                        node (n10) {10}
                    }
                    child {
                        node (n11) {11}
                    }
                }
            }
            child {
                node (n3) {3}
                child {
                    node (n6) {6}
                    child {
                        node (n12) {12}
                    }
                    child {
                        node (n13) {13}
                    }
                }
                child {
                    node (n7) {7}
                    child {
                        node (n14) {14}
                    }
                    child {
                        node (n15) {15}
                    }
                }
            };
        
        \node[box, fit=(n1) (n2) (n3) (n4) (n5) (n6) (n7), draw=blue] {};
        \node[box, fit=(n8) (n9) (n10) (n11), draw=green] {};
        \node[box, fit=(n12) (n13) (n14) (n15), draw=green] {};
\end{tikzpicture}
    }
    \caption{Leaf recomputation strategy}
    \label{fig:naive_recomputation}
\end{subfigure}
\begin{subfigure}[b]{0.5\textwidth}
        
    \centering
    \resizebox{\linewidth}{!}{
    \begin{tikzpicture}[
        level 1/.style={sibling distance=100mm},
        level 2/.style={sibling distance=50mm},
        level 3/.style={sibling distance=25mm},
        every node/.style={circle, draw, minimum size=20pt},
        edge from parent/.style={draw, -},
        box/.style={rectangle, draw=gray, dashed, rounded corners, inner sep=8pt, line width=1.5pt}
    ]
        \node (n1) {1}
            child {
                node (n2) {2}
                child {
                    node (n4) {4}
                    child {
                        node (n8) {8}
                    }
                    child {
                        node (n9) {9}
                    }
                }
                child {
                    node (n5) {5}
                    child {
                        node (n10) {10}
                    }
                    child {
                        node (n11) {11}
                    }
                }
            }
            child {
                node (n3) {3}
                child {
                    node (n6) {6}
                    child {
                        node (n12) {12}
                    }
                    child {
                        node (n13) {13}
                    }
                }
                child {
                    node (n7) {7}
                    child {
                        node (n14) {14}
                    }
                    child {
                        node (n15) {15}
                    }
                }
            };
        
        \node[box, fit=(n1), draw=blue] {};
        \node[box, fit=(n2) (n4) (n5) (n8) (n9) (n10) (n11), draw=green] {};
        \node[box, fit=(n3) (n6) (n7) (n12) (n13) (n14) (n15), draw=green] {};
\end{tikzpicture}
    }
    \caption{Subtree recomputation strategy}
    \label{fig:our_recomputation}
\end{subfigure}
\caption{Comparing the batching patterns of the two different recomputation strategies. The leaf recomputation strategy in \cref{fig:naive_recomputation} performs the local solve stage operations in large batches and then performs all of the merge stages in a separate batch. Our proposed subtree recomputation strategy (\cref{fig:our_recomputation}) performs local solves and multiple levels of the merge stage for a complete subtree of the discretization tree structure.}
\label{fig:graph_view_of_recomputation}
\end{figure}
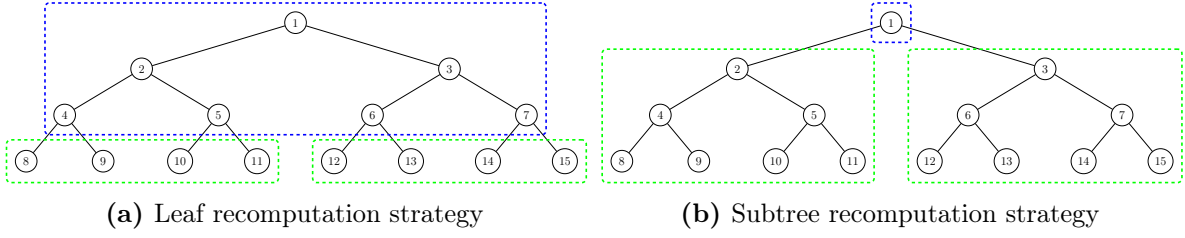

We find that it is advantageous to push the idea of reducing data transfers at the cost of more floating-point operations further.
In our proposed recomputation strategy (\cref{alg:our_recomputation}), we delete and recompute the products of the local solve stage and multiple levels of the merge stage. 
To implement this, we operate in batches defined by ``complete subtrees'', which are subtrees containing all of the descendants of a particular node $j$.
We break the lowest levels of the discretization tree into the largest complete subtrees where the computations in \cref{alg:local_solve,alg:merge} can all fit into a GPU's on-device memory.
The size of these maximal complete subtrees varies depending on GPU memory, polynomial order $p$, and floating-point datatype; in our experiments, these maximal complete subtrees usually have depth $6$ or $7$.
For each such complete subtree, we perform all of the local solve and merge operations, after which we only save the top-level Poincar\'e--Steklov matrix and outgoing boundary data vector. Because this is a small amount of data, we can store it on the GPU, and do not need to move the outputs to host RAM. 
After processing all of the subtrees, the final merge stages are performed on the GPU. 
The downward pass is evaluated sequentially on the different subtrees, at which point the local solve stage and low-level merges must be recomputed.
This recomputation method was inspired by optimizations for contemporary deep learning architectures \citep{dao_flashattention_2022,gu_mamba_2024}, which suggest kernel fusion, a technique that performs multiple steps of a sequential computation at once to keep the necessary data near the processor cores.  
We visualize the different recomputation methods in \cref{fig:graph_view_of_recomputation}.

\begin{figure}[h!]
    \begin{minipage}[b]{0.5\textwidth}
        \centering
    \includeinkscape[width=\linewidth]{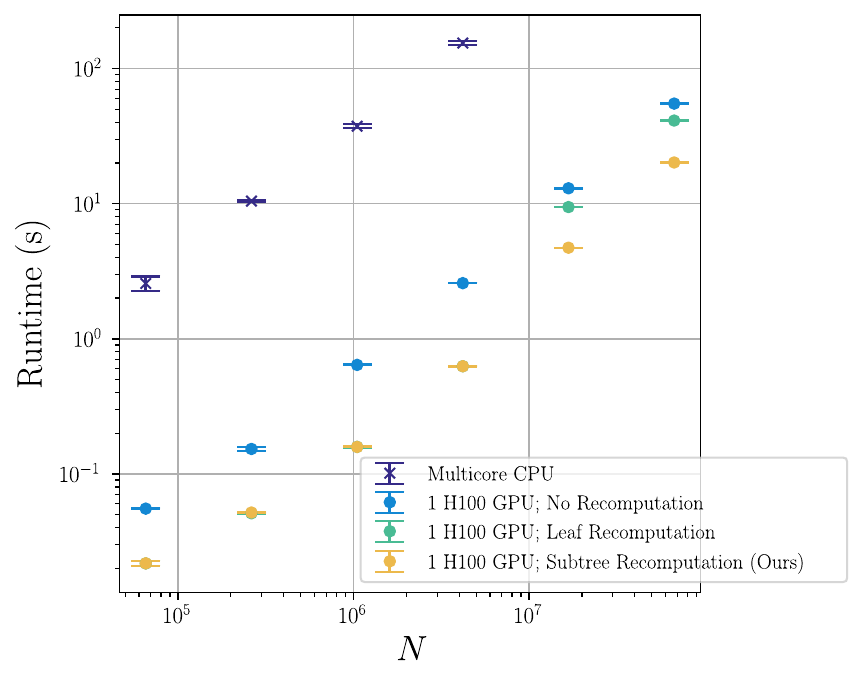}
    \end{minipage}
    \hspace{-3em}
    \begin{minipage}[b]{0.55\textwidth}
        \resizebox{\textwidth}{!}{
            \begin{tabular}{lrrr}
                \textbf{Method} & $N$ & \textbf{Runtime (s)} & \textbf{\% of Peak FLOPS} \\
                \midrule
                No recomputation & $4{,}194{,}304$ & $2.58$ & $1.56\%$ \\
                No recomputation & $16{,}777{,}216$ & $12.99$ & $2.30\%$ \\
                No recomputation & $67{,}108{,}864$ & $54.99$ & $4.17\%$ \\
                \midrule
                Leaf recomputation & $4{,}194{,}304$ & $0.63$ & $6.45\%$ \\
                Leaf recomputation & $16{,}777{,}216$ & $9.43$ & $3.27\%$ \\
                Leaf recomputation & $67{,}108{,}864$ & $41.18$ & $5.66\%$ \\
                \midrule
                Subtree recomputation (Ours) & $4{,}194{,}304$ & $0.63$ & $6.45\%$  \\
                Subtree recomputation (Ours) & $16{,}777{,}216$ & $4.02$ & $14.86\%$  \\
                Subtree recomputation (Ours) & $67{,}108{,}864$ & $17.43$ & $20.01\%$\\
                \bottomrule
            \end{tabular}
            }
            \vspace{6em}
    \end{minipage}
    \caption{Even with a na\"ive implementation of the HPS algorithm which does not perform any recomputation, using a single GPU achieves large speedups over a multicore CPU system. 
    When we use our proposed subtree recomputation strategy, the speedup increases by another factor of two. 
    (Left) \revtwo{We vary $L=4, \dots, 9$  and hold $p=16$ fixed to generate problems with $N=p^2 4^L = 256 \times 4^L = 4^{L+4}$ degrees of freedom.}
    We measure the total runtime of our 2D method merging DtN matrices\revtwo{; this runtime includes the execution of the local solve stage, the merge stage, and the downward pass}.
    Vertical error bars show $\pm$ 1 standard error computed over five trials. 
    (Right) For each of the GPU implementations, we compute the total number of FLOPS and report this as a percentage of the GPU's peak FLOPS, estimated by the manufacturer to be $34 \times 10^{12}$.}
    \label{fig:recomputation_results}
\end{figure}

In \cref{fig:recomputation_results}, we compare the performance of our method across computer architectures and recomputation strategies. 
We consider two different architectures, a multicore Intel Xeon node with a 64-core processor, and a GPU architecture using a single Nvidia H100 GPU.
Evaluating our method on the GPU gives us significant speedups over the multicore CPU architecture, even when using an implementation with no recomputation strategy, which transfers all precomputed matrices to host RAM after each algorithm step. 
The two recomputation strategies begin to diverge for problem sizes over $10^7$ discretization points, at which point the precomputed matrices cannot all fit on the GPU.
\revtwo{We use subtree depth $7$ for this experiment; we explore the effect of this choice in \mbox{\cref{sec:addl_results}} by measuring the runtime for different subtree depths.} 
We also estimate the percentage of peak double-precision floating point operations per second (FLOPS) achieved by the different recomputation strategies.
Our proposed recomputation strategy uses the most floating-point operations and has the fastest runtime, which means it reaches a higher percentage of peak FLOPS than the other implementations.

\begin{algorithm}[h!]
    \caption{
       Leaf recomputation strategy.
    }
    \label{alg:baseline_recomputation}
    \KwIn{Differential operators $\{ \boldsymbol{L}^{(i)} \}_{i=1}^{\nleaves} $; 
    source functions $\{ \boldsymbol{f}^{(i)} \}_{i=1}^{\nleaves}$; 
    boundary data $\boldsymbol{g}$}
    Let $b$ be the maximum batch size that can fit on the GPU \\
    Split the indices $\{ 1, 2, \dots, \nleaves \}$ into batches $I_1, I_2, \dots, I_{\lceil \nleaves / b \rceil}$ \\
    \For(){Batch $j$}{
        Move $\{ \boldsymbol{L}^{(i)} \}_{i \in I_j}$ and  $\{ \boldsymbol{f}^{(i)} \}_{i \in I_j}$ to the GPU \\
        Perform the local solve stage for this batch of leaves \\
        Delete $\{ \boldsymbol{Y}^{(i)} \}_{i \in I_j}$ and  $\{ \boldsymbol{v}^{(i)} \}_{i \in I_j}$ \\
        Transfer $\{ \boldsymbol{T}^{(i)} \}_{i \in I_j}$ and $\{ \boldsymbol{h}^{(i)} \}_{i \in I_j}$ to host RAM
    }
    Concatenate $\{ \boldsymbol{T}^{(i)}\}_{i=1}^{\nleaves}$ and $\{ \boldsymbol{h}^{(i)}\}_{i=1}^{\nleaves}$ and transfer to GPU \\
    Perform all merge operations on the GPU and transfer all $ \boldsymbol{S}^{(i)}$ and $\boldsymbol{\tilde g}^{(i)}$ to host RAM \\
    Propagate boundary data to the leaves \\
    Transfer leaf-level boundary data $\{\boldsymbol{g}^{(i)} \}_{i=1}^{\nleaves}$ to host RAM \\
    \For(){Batch $j$}{
        Move $\{ \boldsymbol{L}^{(i)} \}_{i \in I_j}$, $\{ \boldsymbol{f}^{(i)} \}_{i \in I_j}$, and $\{\boldsymbol{g}^{(i)} \}_{i \in I_j}$ to the GPU \\
        Compute local solutions $\boldsymbol{u}^{(i)}_{i \in I_j}$ \\
        Transfer $\{ \boldsymbol{u}^{(i)} \}_{i \in I_j}$ to host RAM
    }
    \KwResult{Solutions on the Chebyshev discretization points $\{\boldsymbol{u}^{(i)} \}_{i=1}^{\nleaves}$ }
\end{algorithm}

\begin{algorithm}[h!]
    \caption{
       Subtree recomputation strategy.
    }
    \label{alg:our_recomputation}
    \KwIn{Differential operators $\{ \boldsymbol{L}^{(i)} \}_{i=1}^{\nleaves} $; 
    source functions $\{ \boldsymbol{f}^{(i)} \}_{i=1}^{\nleaves}$; 
    boundary data $\boldsymbol{g}$}
    Let $M = {m_1, m_2, ...}$ be the roots of the maximal subtrees \\
    \For(){Subtree rooted at  $m_j$}{
        Let $I_j$ be the set of leaves of the subtree \\
        Move $\{ \boldsymbol{L}^{(i)} \}_{i \in I_j}$ and  $\{ \boldsymbol{f}^{(i)} \}_{i \in I_j}$ to the GPU \\
        Perform the local solve stage for this subtree \\
        Delete $\{ \boldsymbol{Y}^{(i)} \}_{i \in I_j}$ and  $\{ \boldsymbol{v}^{(i)} \}_{i \in I_j}$ \\
        Merge the leaves to the top of subtree $m_j$, deleting all outputs except $\boldsymbol{T}^{(m_j)}$ and $\boldsymbol{h}^{(m_j)}$ \\
        Keep $\boldsymbol{T}^{(m_j)}$ and $\boldsymbol{h}^{(m_j)}$ on GPU \\
    }
    Perform final merge operations on the GPU \\
    Propagate boundary data to the roots of the maximal subtrees \\
    Transfer boundary data to host RAM \\
    \For(){Subtree rooted at  $m_j$}{
        Let $I_j$ be the set of leaves of the subtree \\
        Move $\{ \boldsymbol{L}^{(i)} \}_{i \in I_j}$, $\{ \boldsymbol{f}^{(i)} \}_{i \in I_j}$ and $\boldsymbol{g}^{(m_j)}$ to the GPU \\
        Perform the local solve stage for this subtree \\
        Merge the leaves to the top of subtree $j$ \\
        Propagate boundary information down the tree to the leaves \\
        $\boldsymbol{u}^{(i)} \gets \boldsymbol{Y}^{(i)} \boldsymbol{g}^{(i)} + \boldsymbol{v}^{(i)} $ \\
        Transfer $\{ \boldsymbol{u}^{(i)} \}_{i\in I_j}$ to host RAM
    }
    \KwResult{Solutions on the Chebyshev discretization points $\{\boldsymbol{u}^{(i)} \}_{i=1}^{\nleaves}$ }
\end{algorithm}

\subsection{An adaptive discretization strategy to reduce memory complexity in 3D}
\label{sec:hardware_acceleration_3D}
When extending from two to three dimensions, we face different computational challenges. 
A large part of the difficulty involves the size of the matrices arising in the merge stage discussed in \cref{sec:hps_algos_merge}.
For each merge operation, the matrix $\boldsymbol{D}$ must be inverted. 
This matrix has a number of rows and columns proportional to the number of discretization points that lie along the interfaces being merged. 
\revone{In two dimensions, the size of this merge interface is $O\left( p2^\ell \right)$ to merge nodes $\ell$ levels above the leaves. In three dimensions, the size of this merge interface is $O\left(p^2 4^\ell \right)$.}
Because this quantity grows very quickly \revone{for 3D problems} as we increase the tree depth $L$, we quickly run out of memory required to store and invert the matrices on the GPU at the highest level of the merge stage. 
Performing the top-level merge operation is when the instantaneous memory footprint peaks, as space for $\boldsymbol{D}$, $\boldsymbol{D}^{-1}$, and various buffers must all be allocated on the GPU simultaneously. In contrast with other parts of the algorithm, this peak memory footprint is not reducible by strategies such as batching or data transfer. 
\revone{Because of this need to invert matrices near the memory limits of the GPU, we transfer data to and from the GPU at each merge level and do not use the recomputation strategies discussed in \mbox{\cref{sec:hardware_acceleration_2D}}.}

To reduce the size of $\boldsymbol{D}$ at the final merge step, we propose to extend the adaptive HPS method presented for 2D problems in \cite{geldermans_adaptive_2019} to three dimensions. 
This method adaptively refines element sizes in a data-dependent manner, in an effort to concentrate discretization points in the regions of the domain where the coefficient and source functions have high local variation. 
In \cref{tab:merge_sizes}, we show that the size of $\boldsymbol{D}$ generated using our adaptive refinement technique is much smaller than that of the uniform refinement with no loss in accuracy.

Developing a version of the HPS methods presented in \cref{sec:hps_algos} that is compatible with an adaptive discretization requires slight modification of the algorithms presented in the previous section. The major changes are the introduction of a method for adaptively refining our octree and a method for merging nodes with different levels of refinement.

\begin{table}[!ht]
    \centering
    \begin{tabular}{
        lrrrr
        }
        \toprule
            \textbf{Method} & $\boldsymbol{p}$ & \textbf{Relative} $\boldsymbol{\ell_\infty}$ \textbf{Error} & \textbf{\# Leaves} & \textbf{Size of} $\boldsymbol{D}$    \\
        \midrule
        Uniform & 8 & $1.48\times 10^{-4}$ & $512$  & $6{,}912$\\
        Adaptive & 8 & $1.45\times 10^{-4}$  & $190$ & $2{,}700$ \\ %
        \midrule
        Uniform & 12 & $3.62 \times 10^{-7}$ & $512$ & $19{,}200$ \\
        Adaptive & 12 & $2.04 \times 10^{-7}$ & $442$  & $7{,}500$ \\ %
        \midrule
        Uniform & 16 & $4.20 \times 10^{-6}$ & $64$ & $9{,}408$ \\
        Adaptive & 16 & $1.41 \times 10^{-6}$ & $57$ & $4{,}116$ \\ %
        \bottomrule
    \end{tabular}
    \caption{Adaptive discretization methods can greatly reduce the peak memory requirements of HPS methods in three dimensions. 
    We present the size of the discretization tree and final merge steps in our 3D ``wavefront'' example (\cref{sec:3D_wave_front}).
    We compare adaptive and uniform discretizations that have similar errors and observe the adaptive discretization strategy can greatly reduce the size of the final $\boldsymbol{D}$ matrix, a proxy for peak memory usage. 
    }
    \label{tab:merge_sizes}
\end{table}

\subsubsection{Criterion for adaptive refinement}

For a given leaf of the discretization tree, let $\bx_0$ be the set of $p^3$ Chebyshev points discretizing the leaf. Let $\bx_1$ be the set of $8p^3$ discretization points found by breaking the leaf into eight children and creating a Chebyshev grid on each child. 
Let $\boldsymbol{L}_{8f1}$ be an interpolation matrix mapping from $\bx_0$ to $\bx_1$. We evaluate whether a function is sufficiently refined on a leaf by checking whether we can use polynomial interpolation to accurately map from evaluations on $\bx_0$ to evaluations on $\bx_1$, relative to the global $L_\infty$ norm of the function.
We specify a tolerance parameter $\epsilon$, and for each leaf in our tree, we check the following condition:
\begin{align}
    \frac{ \| f(\bx_1) - \boldsymbol{L}_{8f1} f(\bx_0) \|_\infty}{ \| f \|_\infty} < \epsilon.
    \label{eq:refinement_check}
\end{align}
If this condition is met, we say the leaf is sufficiently refined. 
Otherwise, we split the leaf into eight children and check each child. 
We form a final discretization tree by refining each coefficient function in our differential operator, as well as the source term, and taking the union of the resulting trees.
Additionally, we refine a few extra leaves to enforce a ``level restriction'' criterion, which specifies that no leaf can have a side length greater than twice that of its neighbors. This method is an extension of the method for two-dimensional problems presented in \cite{geldermans_adaptive_2019}, which uses a similar relative $L_2$ convergence criterion and level restriction criterion.

\subsubsection{Local solve stage}
The local solve stage for adaptively refined discretization trees is the same as in the uniform refinement case. 
Although they are defined over leaves with different volumes, each local boundary value problem has the same number of interior and boundary discretization points, and the local problems are still embarrassingly parallel. Thus, we can use batched linear algebra to accelerate this part of the algorithm. 

\subsubsection{Merging nodes with different discretization levels}
The nonuniform merge stage is different from the uniform merge stage because neighboring nodes may have different refinement levels, which means the discretization points along either side of the merge interface may not exactly align. 
Recall the block linear system (\cref{eq:merge_lin_system}) arising during the merge stage:
\begin{linenomath}\begin{align*}
    \begin{bmatrix} \boldsymbol{A} \vphantom{\boldsymbol{g}_{\ext}^{(j)}} & \boldsymbol{B} \\
        \boldsymbol{C} & \boldsymbol{D}
    \end{bmatrix} 
    \begin{bmatrix}
        \boldsymbol{g}_{\ext}^{(j)} \\ \boldsymbol{g}_{\intt}^{(j)}
    \end{bmatrix} = \begin{bmatrix}
        \boldsymbol{u}^{(j)}_{\ext} -\boldsymbol{h}^{(\child)}_{\ext} \\
        - \boldsymbol{h}_{\intt}^{(\child)}
    \end{bmatrix}.
\end{align*}\end{linenomath}
When neighboring volume elements have different refinement levels, there will be a mismatch between the interior boundary discretization points on either side of the merge interface. 
We need to decide how to represent $\boldsymbol{g}_{\intt}^{(j)}$, $\boldsymbol{h}_{\intt}^{(\child)}$, $\boldsymbol{B} $, $\boldsymbol{C}$, and $\boldsymbol{D}$ in  \cref{eq:merge_lin_system}. 
We choose to discretize these objects using the coarser of the two sets of discretization points along the merge interface; the discretization points along the exterior boundary elements are inherited from the child nodes. 
To assemble the blocks in \cref{eq:merge_lin_system}, this requires projecting some rows and columns of the Poincar\'e--Steklov matrices using precomputed interpolation operators which map between one and four 2D Gauss--Legendre panels. The ``level restriction''  constraint greatly simplifies this compression because the resulting projection operations are guaranteed to be four-to-one.

\subsubsection{Downward pass}
To propagate the boundary information to the leaf nodes, we follow the general structure of \cref{alg:down_pass}. However, we must undo the projection along merge interfaces that occurs during the nonuniform merge stage. This is accomplished by applying the precomputed interpolation operators to the boundary data $\boldsymbol{g}^{(i)}$.

\section{Numerical examples in two dimensions}
\label{sec:2D_examples}
In this section, we present numerical results on problems with two spatial dimensions. 
All experiments in this section were conducted using one Nvidia H100 GPU and a host memory space with 100GB of RAM. 
In all of the experiments, we use the novel subtree recomputation strategy introduced in \cref{sec:hardware_acceleration_2D}\revtwo{; the DtN version of this recomputation strategy uses subtrees of depth 7, and the ItI version of this recomputation strategy uses subtrees of depth 6.}

\subsection{High-order convergence on variable-coefficient problems with known solutions}
\begin{figure}[!ht]
    \centering
    \begin{subfigure}[b]{0.5\textwidth}
        \centering
        \includeinkscape[height=6cm]{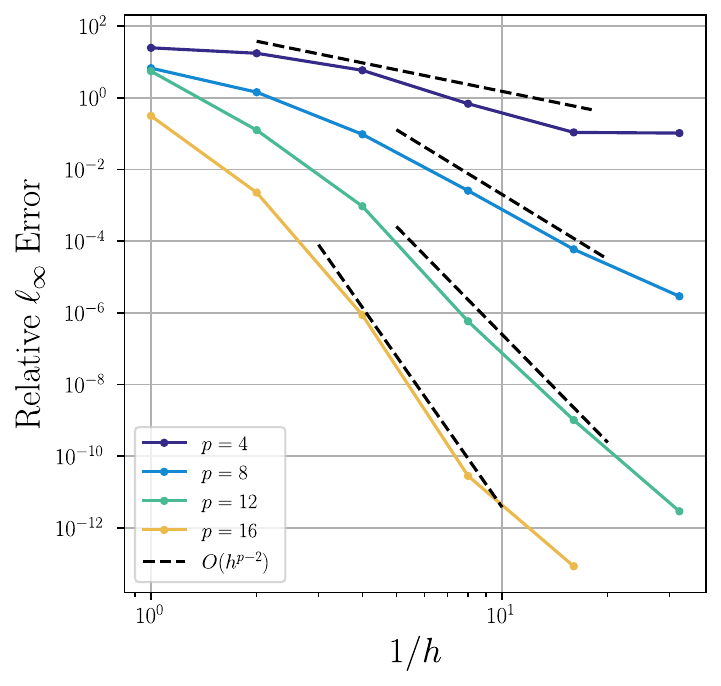}
        \caption{DtN}
        \label{fig:hp_convergence_DtN}
    \end{subfigure}%
    \begin{subfigure}[b]{0.5\textwidth}
        \centering
        \includeinkscape[height=6cm]{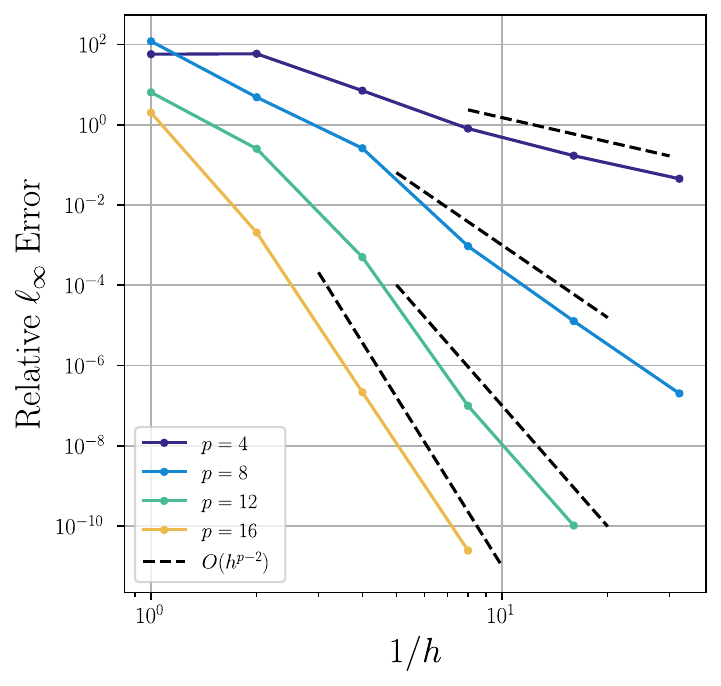}
        \caption{ItI}
        \label{fig:hp_convergence_ItI}
    \end{subfigure}
    \caption{Using $p$ Chebyshev points per dimension on each leaf, and leaves of size $h$, the relative $\ell_\infty$ errors of our method converge at rate $O(h^{p-2})$. \cref{fig:hp_convergence_DtN} shows the convergence of the HPS method using DtN matrices applied to \cref{eq:DtN_problem} and \cref{fig:hp_convergence_ItI} shows the convergence of the HPS method using ItI matrices applied to \cref{eq:ItI_problem}.}
    \label{fig:hp_convergence}
\end{figure}

There are two main ways to increase the accuracy of our composite spectral colocation scheme: refine each leaf patch into four children or increase the polynomial order of the representation of the solution on each patch. 
Empirically, the error is controlled by the polynomial order $p$ and the side length of each leaf $h$.
In our implementation, $p$ is specified by the user and $h$ is controlled by $L$, the user-specified depth of the discretization tree. 
The tradeoff between these two parameters is a widely-studied topic in numerical analysis and goes by the name of ``$hp$-adaptivity''. 
We study the $hp$-adaptivity properties of our solver using two problems with variable-coefficient differential operators and known solutions. The first problem is a variant of Poisson's equation with spatially-varying coefficients:
\begin{align}
\left\{ 
    \begin{aligned}
         \Delta u(x) - \cos(5 x_2)u_{x_1}(x) + \sin(5 x_2) u_{x_2}(x) &= f(x), &&\quad x \in [-1,1]^2, \\
        u(x) &= g(x), &&\quad x \in \partial [-1,1]^2,
    \end{aligned}
\right.
    \label{eq:DtN_problem}
\end{align}
where $u_{x_1}$ and $u_{x_2}$ are the partial derivatives of $u$ in directions $x_1$ and $x_2$, respectively. 
We manufacture the source $f$ and the Dirichlet data $g$ so the solution to this problem is 
\begin{linenomath}\begin{align*}
    u(x) = u(x_1, x_2) = e^{5x_1} \sin(5x_2) + \sin(10 \pi x_1) \sin(\pi x_2).
\end{align*}\end{linenomath}
We also study an inhomogeneous Helmholtz problem with a Robin boundary condition:
\begin{align}
\left\{
    \begin{aligned}
         \Delta u(x) + (1 + e^{-50\| x \|^2}) u(x)  &= f(x), && \quad x \in [-1,1]^2, \\
         u_n(x)  + i u(x) &= g(x), && \quad x \in \partial [-1,1]^2.
    \end{aligned}
    \right.
    \label{eq:ItI_problem}
\end{align}
We manufacture the source $f$ and the Robin data $g$ so solution to this problem is 
\begin{linenomath}\begin{align*}
    u(x) = u(x_1,x_2) = e^{i 20 x_1} + e^{i30x_2}.
\end{align*}\end{linenomath}
\cref{fig:hp_convergence} shows the convergence of our method on these problems. 
We measure the relative error of our computed solution $\boldsymbol{u}$ by computing $\| \boldsymbol{u} - \boldsymbol{u}_{\text{true}} \|_\infty / \| \boldsymbol{u}_{\text{true}} \|_\infty$. 
The $\ell_\infty$ norms are estimated by taking the maximum over all interior discretization points. 
In \cref{fig:hp_convergence}, we see the errors of both the DtN and ItI versions of the method converging at rate $O(h^{p-2})$, even for high polynomial orders.

\subsection{High-frequency forward wave scattering problems}
In this example, we solve a variable-coefficient Helmholtz equation coupled with a Sommerfeld radiation condition. The system is excited by a plane wave with direction $\hat{s} = [1, 0]^\top$ and frequency $k$:
\begin{equation}
\left\{
    \begin{aligned}
        \Delta u(x) + k^2(1 + q(x)) u(x) &= -k^2 q(x) e^{ik\langle \hat{s}, x \rangle}, && \quad x \in [-1, 1]^2, \\
        \sqrt{r} \left( \tfrac{\partial u}{\partial r} - iku \right) &\to 0, && \quad r = \|x \|_2 \to \infty.
    \end{aligned}
\right.
    \label{eq:forward_scattering_problem}
\end{equation}
\cref{eq:forward_scattering_problem} models time-harmonic wave scattering in many imaging modalities, such as sonar or radar imaging, geophysical sensing, and nondestructive testing of materials \citep{borges_high_2016}. 
As such, forward wave scattering solvers are often used inside an inner loop of optimization routines for solving these inverse problems. 
These forward solvers must be highly optimized as they are evaluated hundreds or thousands of times over the course of an algorithm.

To solve \cref{eq:forward_scattering_problem}, we use the ItI variant of our implementation, and at the top level of the merge stage, we solve a boundary integral equation which enforces the Sommerfeld radiation condition \citep{gillman_spectrally_2015}. 
Discretizing the boundary integral equation requires a high-order Nystr\"om method to generate single- and double-layer potential matrices, which we perform in MATLAB using the \texttt{chunkIE} package \citep{askham_chunkie_2024}. 
For each discretization level and frequency, generating these matrices takes a few seconds on a standard laptop. 
Because these matrices can be precomputed once for a given discretization level and frequency, we do not include the time required to generate these matrices in our runtime measurements. 
The solution of this boundary integral equation specifies incoming impedance data, which is propagated down the tree to form the interior solution. 
While we solve \cref{eq:forward_scattering_problem} for one source direction $\hat s$, this scheme can compute solutions for multiple different sources in parallel at the cost of a few extra matrix-vector multiplications.

\begin{figure}[!ht]
    \begin{subfigure}[b]{0.33\textwidth}
        \centering
        \includeinkscape[height=4cm]{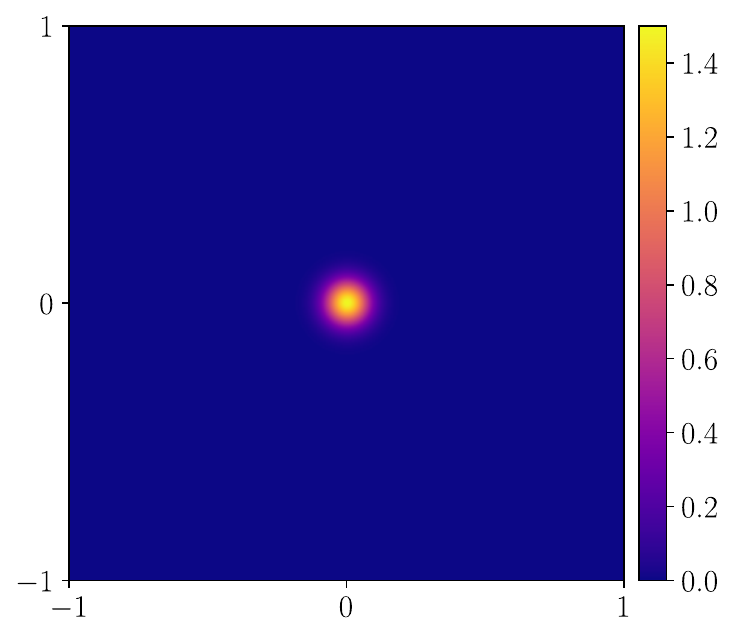}
        \caption{Scattering potential $q(x)$}
        \label{fig:q_GBM_1}
    \end{subfigure}%
    \begin{subfigure}[b]{0.33\textwidth}
        \centering
        \includeinkscape[height=4cm]{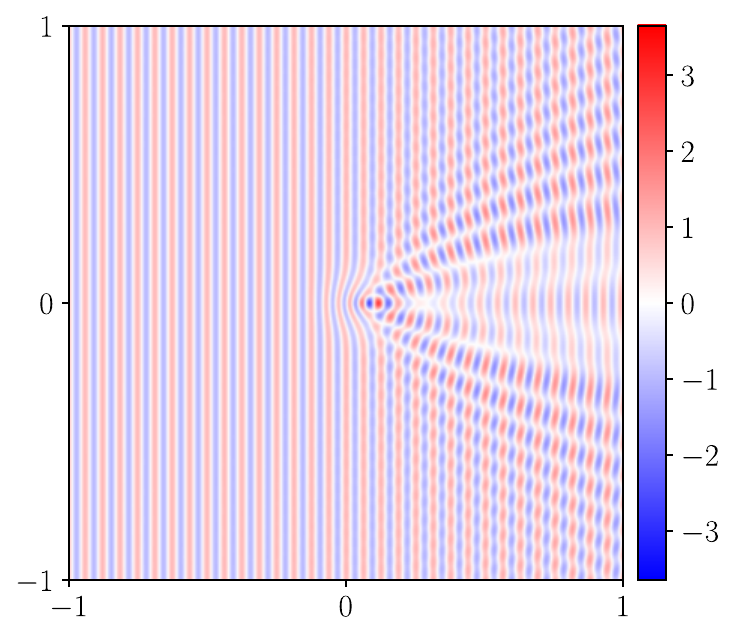}
        \caption{Total wave field; $k=100$}
        \label{fig:k100_utot_GBM_1}
    \end{subfigure}
    \begin{subfigure}[b]{0.33\textwidth}
        \centering
        \includeinkscape[height=4cm]{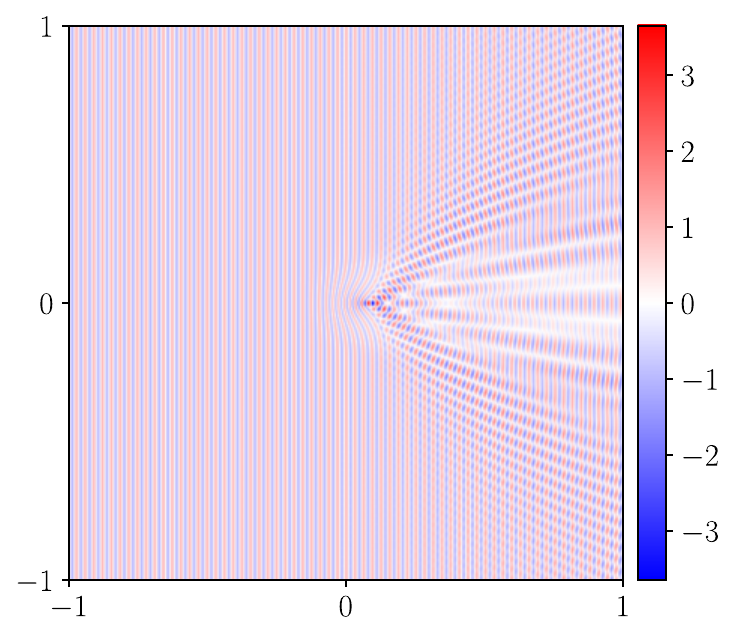}
        \caption{Total wave field; $k=200$}
        \label{fig:k200_utot_GBM_1}
    \end{subfigure}
    \caption{Visualizing the solutions of forward scattering problems for the single Gaussian bump scattering potential. \cref{fig:k100_utot_GBM_1,fig:k200_utot_GBM_1} show the real part of the total wave field.}
    \label{fig:forward_scattering_examples_GBM_1}
\end{figure}

\begin{figure}[!ht]
    \centering
    \begin{subfigure}[b]{0.5\textwidth}
        \centering
        \includeinkscape[height=6cm]{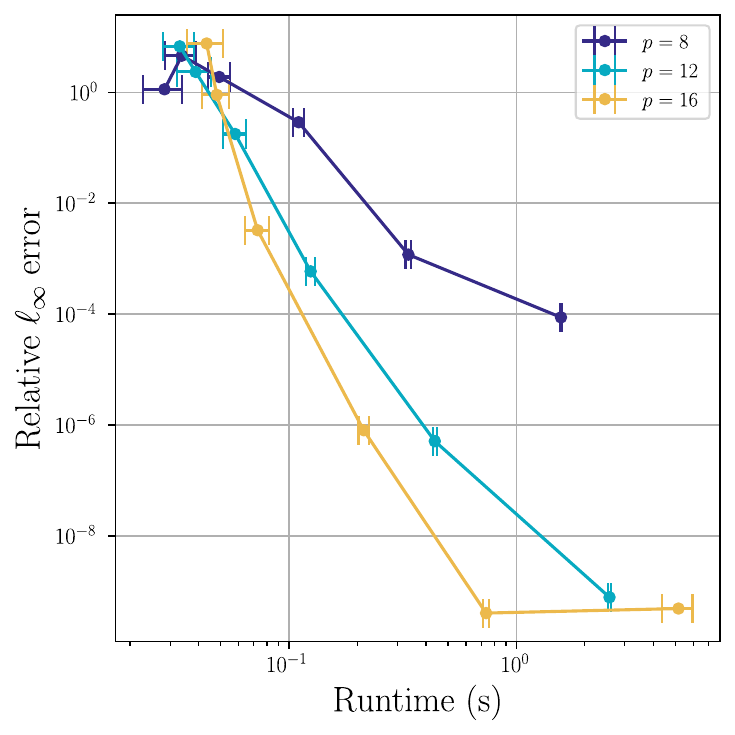}
        \caption{Runtime and errors when $k=100$}
        \label{fig:k100_scattering_GBM_1}
    \end{subfigure}%
    \begin{subfigure}[b]{0.5\textwidth}
        \centering
        \includeinkscape[height=6cm]{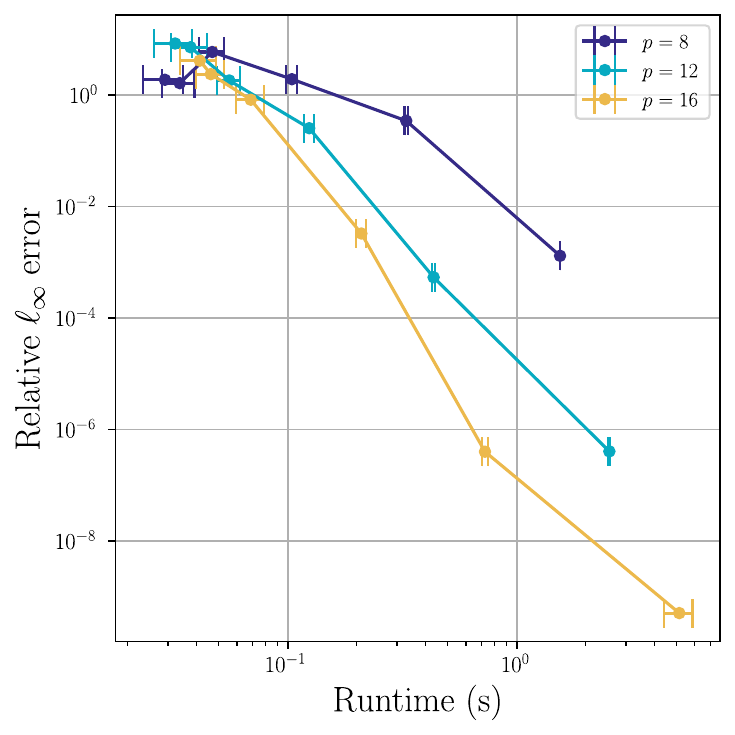}
        \caption{Runtime and errors when $k=200$}
        \label{fig:k200_scattering_GBM_1}
    \end{subfigure}
    \caption{Error-runtime study on the single Gaussian bump scattering potential. Using GPU acceleration, our method can rapidly converge to high-accuracy solutions in high-frequency wave scattering problems. 
    The runtime measurements include the runtime of the entire HPS algorithm and the setup and solution of the boundary integral equation enforcing the radiation condition. Horizontal error bars show $\pm 1$ standard error computed over five trials.
    }
    \label{fig:forward_scattering_convergence_GBM_1}
\end{figure}

\begin{figure}[!ht]
    \begin{subfigure}[b]{0.33\textwidth}
        \centering
        \includeinkscape[height=4cm]{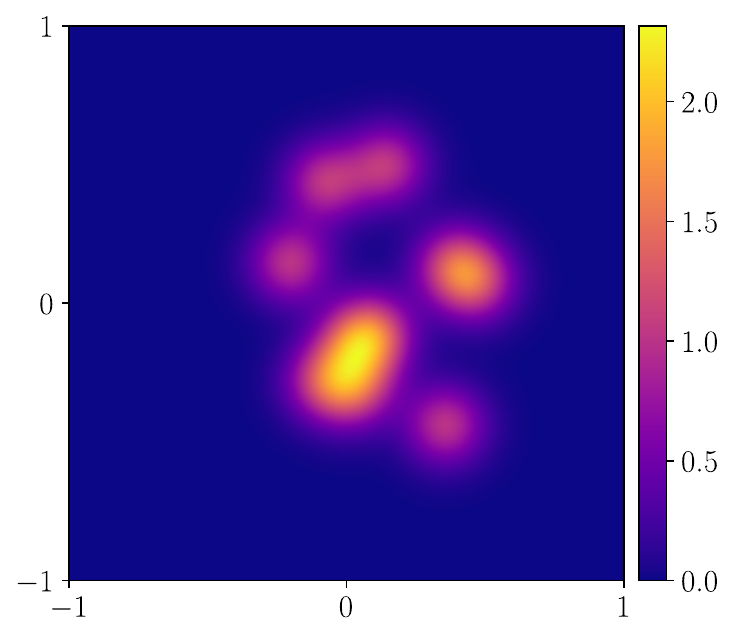}
        \caption{Scattering potential $q(x)$}
        \label{fig:q}
    \end{subfigure}%
    \begin{subfigure}[b]{0.33\textwidth}
        \centering
        \includeinkscape[height=4cm]{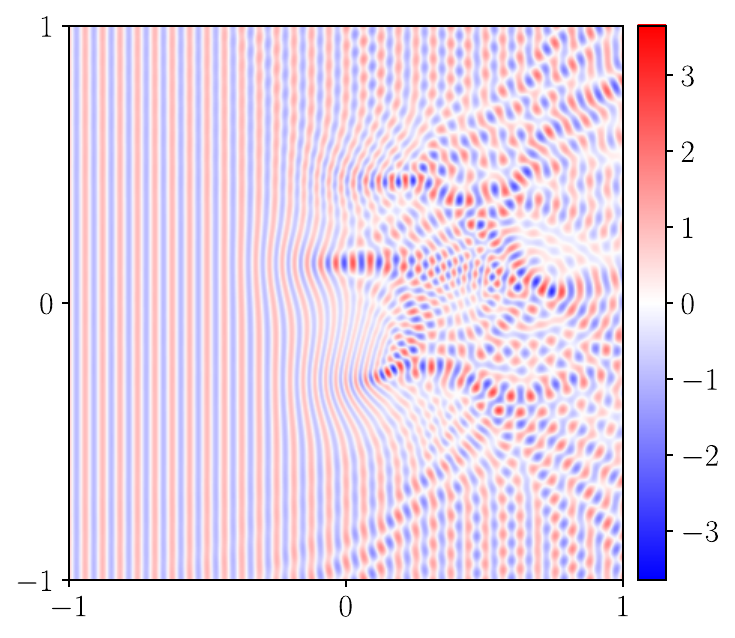}
        \caption{Total wave field; $k=100$}
        \label{fig:k100_utot}
    \end{subfigure}
    \begin{subfigure}[b]{0.33\textwidth}
        \centering
        \includeinkscape[height=4cm]{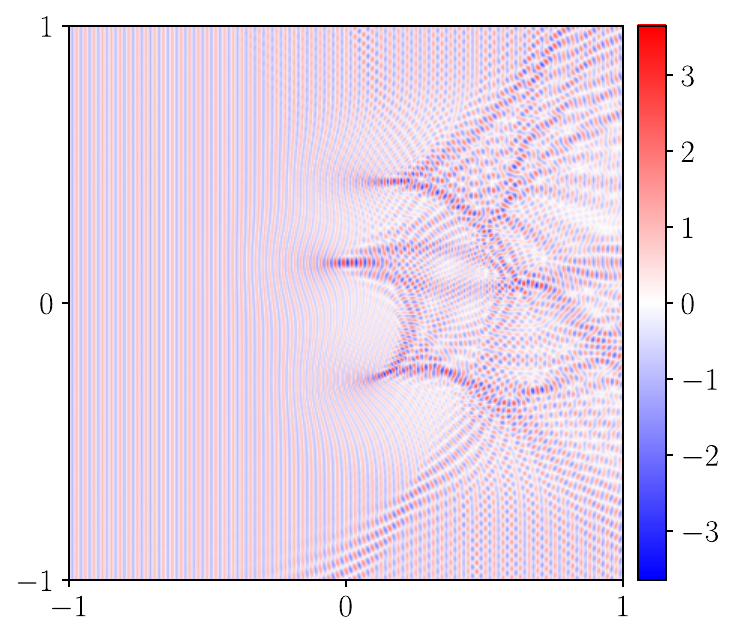}
        \caption{Total wave field; $k=200$}
        \label{fig:k200_utot}
    \end{subfigure}
    \caption{Visualizing the solutions of forward scattering problems for the sum of randomly placed Gaussian bumps scattering potential. \cref{fig:k100_utot,fig:k200_utot} show the real part of the total wave field.}
    \label{fig:forward_scattering_examples}
\end{figure}

\begin{figure}[!ht]
    \begin{subfigure}[b]{0.5\textwidth}
        \centering
        \includeinkscape[height=6cm]{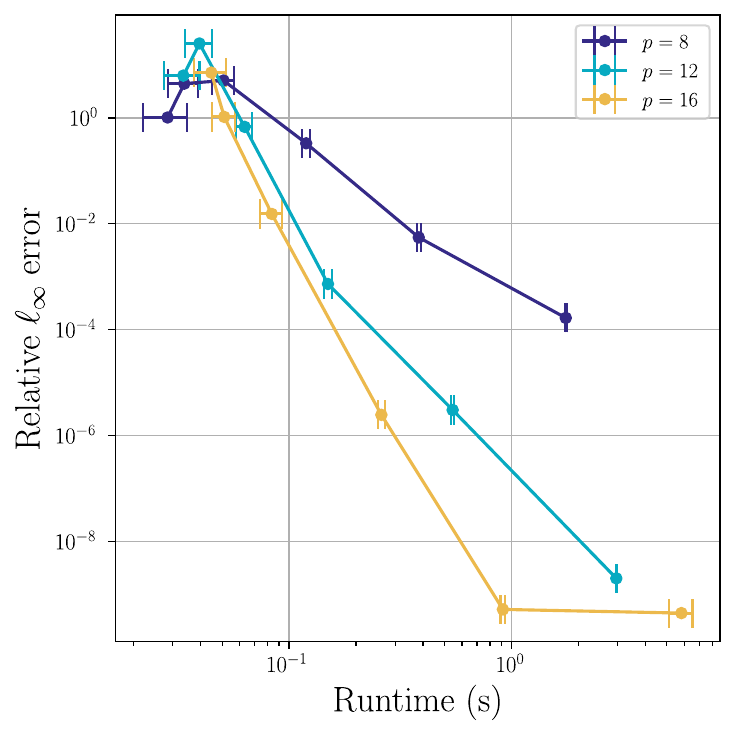}
        \caption{Runtime and errors when $k=100$}
        \label{fig:k100_scattering}
    \end{subfigure}%
    \begin{subfigure}[b]{0.5\textwidth}
        \centering
        \includeinkscape[height=6cm]{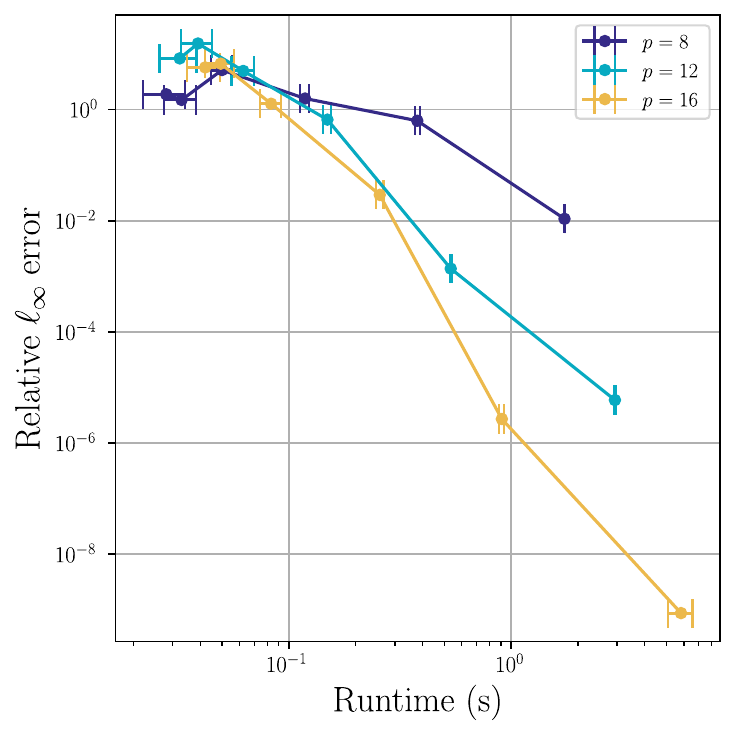}
        \caption{Runtime and errors when $k=200$}
        \label{fig:k20  0_scattering}
    \end{subfigure}
    \caption{Error-runtime study on the sum of randomly placed Gaussian bumps scattering potential. 
    }
    \label{fig:forward_scattering_convergence}
\end{figure}

\renewcommand{\over}{\textrm{over}}
In \cref{fig:forward_scattering_convergence_GBM_1,fig:forward_scattering_convergence}, we measure the runtime and accuracy of our GPU-accelerated solver. Because analytical solutions for \cref{eq:forward_scattering_problem} are unavailable for general scattering potentials $q(x)$, we measure error relative to an overrefined reference solution $\boldsymbol{u}_{\over}$ with approximately $2{,}800$ discretization points in both dimensions.
For each computed solution $\boldsymbol{u}$, we compute the relative $\ell_\infty$ error $\| \boldsymbol{u} - \boldsymbol{u}_{\over} \|_\infty / \| \boldsymbol{u}_{\over} \|_\infty$. 
The $\ell_\infty$ norm is estimated by taking the maximum over a grid of $500\times 500$ regularly-spaced grid points. 
We repeat this experiment for two different choices of the scattering potential $q(x)$. We first choose a single Gaussian bump, which loosely focuses the incoming wave:
\begin{linenomath}\begin{align*}
    q(x) &= 1.5 e^{-160 \| x \|^2 }.
\end{align*}\end{linenomath}
We also consider a collection of randomly placed Gaussian bumps with centers $\{ z^{(i)} \}_{i=1}^{10}$, which causes multiple scattering effects at high wavenumbers:
\begin{align}
    q(x) &= \sum_{i=1}^{10} e^{-50 \| x - z^{(i)} \|^2}. \label{eq:q_sum_gauss}
\end{align}
In \cref{fig:forward_scattering_examples,fig:forward_scattering_examples_GBM_1}, we show these scattering potentials as well as the resulting total wave field $u(x) + e^{ik\langle \hat{s}, x \rangle}$ for different choices of $k$.

Our GPU-accelerated method with our proposed recomputation strategy is very fast. 
It is able to compute high-accuracy solutions to these challenging scattering problems in a few seconds. 
We note that because this method does not rely on any iterative algorithms, the runtime for a fixed discretization level does not depend on the frequency $k$ or the scattering potential $q(x)$. 
However, finer discretizations are required to achieve a fixed error tolerance as $k$ increases. \cref{fig:forward_scattering_convergence_GBM_1,fig:forward_scattering_convergence} show that polynomial order is $p=16$ dominates the lower-order methods for all values of $k$ and scattering potentials considered.
\revtwo{In these plots, we vary the number of degrees of freedom $N=p^2 4^L$ by varying both $p$ and $L$. Different colors correspond to different polynomial orders $p \in \{ 8, 12, 16 \}$, and for a given polynomial order, we vary $L \in \{ 2, 3, \dots, 7 \}$. This generates problems with degrees of freedom varying between $512$ and $4{,}194{,}304$.}

\subsection{Solving inverse scattering problems with automatic differentiation}
\newcommand{\cB}{\mathcal{B}}
Our solver is compatible with the JAX automatic differentiation framework, which allows us to very easily implement gradient-based optimization algorithms using our accelerated HPS solver as a forward model. 
We consider an inverse scattering task to recover \revone{an inhomogeneous scattering potential $q_\theta$ specified by basis coefficients $\{ \theta_j \}$}:
\begin{align}
    \revone{q_\theta(x)} & \revone{= \sum_{b_j \in B_\gamma }\theta_j b_j(x),}
    \label{eq:q_theta} \\
    \revone{B_\gamma} & \revone{= \bigg\{  \sin \left( m_1 \frac{ \pi}{2} \left( x_1  - 1 \right) \right) \sin \left( m_2 \frac{ \pi}{2} \left( x_2  - 1 \right) \right) \bigg| \sqrt{m_1^2 + m_2^2} \leq \gamma \bigg \}.} \label{eq:basis_def}
\end{align}
\revone{The basis $B_\gamma$ is used in \mbox{\citet{borges_high_2016}} because it spans smooth, bandlimited functions which vanish on the boundary of $\Omega$.}
The goal of this inverse scattering task is to estimate $\theta$. \cref{eq:forward_scattering_problem} describes how $q_\theta(x)$ affects the scattered wave field $u_\theta(x)$.
\revone{Our forward model $\mathcal{F} \circ \mathcal{B}$ is a composition of maps where $\cB: \theta \mapsto q_\theta$ is a discrete sine transform and $\cF : q_\theta \mapsto u_\theta$ evaluates the solution of \mbox{\cref{eq:forward_scattering_problem}} at points $x^{(j)}$ away from the support of the scattering potential.}
We choose a ground-truth \revone{$\theta^*$ to be the basis coefficients of the sum of Gaussian bumps scattering potential \mbox{\cref{eq:q_sum_gauss}} projected onto $B_\gamma$} and generate data 
$\mathcal{F}[\theta^*]$. Given this data and our knowledge of the forward model, we wish to recover an estimate of $\theta^*$. We can phrase this problem in an optimization framework: 
\begin{align}
    \argmin_{\theta \in \revone{ \R^{N_\theta}}} \| \revone{\left( \cF \circ \cB \right) (\theta) -  \left( \cF \circ \cB \right) (\theta^*) }  \|^2_2.
    \label{eq:inverse_problem_objective}
\end{align}
We evaluate $\mathcal{F}$ using our GPU-accelerated fast direct solver. Because our solver is compatible with automatic differentiation, we can solve this optimization problem using gradient-based methods. 
\revone{We first evaluate the accuracy of JAX automatic differentiation by comparing the outputs with the action of $J[\theta]$, the Fr\'echet derivative of $\cF \circ \cB$ centered at $\theta$. \mbox{\citet{borges_high_2016}} characterize this derivative in terms of the solution of linear elliptic PDEs which we can compute with our HPS method; we re-state these results in \mbox{\cref{appendix:actions}} for completeness. 
\mbox{\cref{fig:autodiff_accuracy}} shows the convergence of the outputs of automatic differentiation to the Fr\'echet derivative when holding polynomial order $p=16$ fixed and varying the depth of the quadtree, $L=1,\hdots, 5$. 
For leaf size $h$, we observe high-order convergence at rate $O\left(h^{p-2}\right)$. For brevity, we defer the details of this experiment to \mbox{\cref{app:adjoints_experimental_details}}.
We remark that high-order accuracy is possible because $B_\gamma$ defines a smooth global basis which can be represented using the HPS discretization. 
This is in contrast to differentiating with respect to the values of $q$ at the HPS discretization points, which effectively introduces nonsmooth coefficient functions.
We also note that the use of standard automatic differentiation software requires the presence of the entire computational graph in memory, which means it can not be used in conjunction with the recomputation strategies introduced in \mbox{\cref{sec:hardware_acceleration_2D}}.
}

\revone{To solve \mbox{\cref{eq:inverse_problem_objective}}, we propose to} use Gauss--Newton iterations for nonlinear least squares problems. 
This algorithm requires access to $J[\theta]$, \revone{which we implement using JAX automatic differentiation applied to our HPS solver. We use automatic differentiation to implement the actions $J[\theta]^*f$ and $J[\theta]v$ for arbitrary vectors $f,v$ and estimates $\theta$.}
We pair these subroutines with a sparse linear algebra least-squares solver \citep{paige_lsqr_1982} from the SciPy library \citep{2020SciPy-NMeth} to implement the Gauss--Newton algorithm presented in \cref{alg:gauss_newton}.
\begin{algorithm}[h!]
    \caption{
       Gauss--Newton iterations for nonlinear least squares problems
    }
    \label{alg:gauss_newton}
    \KwIn{Data \revone{$\left( \cF \circ \cB \right) (\theta^*)$}; Initial estimate $\theta_0$}
    $t \gets 0$ \\
    \While(){not converged} {
        Compose automatic differentiation and our fast direct solver to compile the function \revone{$f \mapsto  J[\theta_t]^* f$} \\
        Compose automatic differentiation and our fast direct solver to compile the function \revone{$v \mapsto  J[\theta_t]v$} \\
        Define a linear operator $J[\theta_t]$ using subroutines  \revone{$J[\theta_t]^*f$} and $J[\theta_t]v$ \\
        Use a least-squares solver to compute $\delta \gets \argmin\limits_{\delta} \| \revone{ \left( \cF \circ \cB \right) (\theta^*)} - \left( \revone{\left( \cF \circ \cB \right) (\theta_t)} +  J[\theta_t] \delta  \right) \|^2_2$  \\
        $\theta_{t+1} \gets \theta_t + \delta$ \\
        $t \gets t+ 1$ \\
    }
    \KwResult {Final estimate $\theta_t$}
\end{algorithm}

\revone{To solve the optimization problem, we use $\mathcal{B}_\gamma$, $\gamma = 5$, which gives us $N_\theta = 15$ optimization variables. Following \mbox{\citet{borges_high_2016}}, we initialize the optimization variables corresponding to the lowest three frequency components to the ground-truth $\theta_0 = \theta_*$, and we initialize the other variables at $\boldsymbol{0}$. We run 21} iterations of the Gauss--Newton algorithm. 
\cref{fig:inverse_scattering} shows the optimization variables take some iterations to approach $\theta^*$ and then converge superlinearly to \revone{near} machine precision.
Calculating each Gauss--Newton update is fast because of the GPU acceleration of the forward model, \revone{$J[\theta_t]^*f$, and $J[\theta_t]v$}. The entire experiment runs in $75$ seconds using one H100 GPU.

\begin{figure}[!ht]
    \centering
    \begin{subfigure}[b]{0.3\textwidth}
        \centering
        \includegraphics[height=4cm]{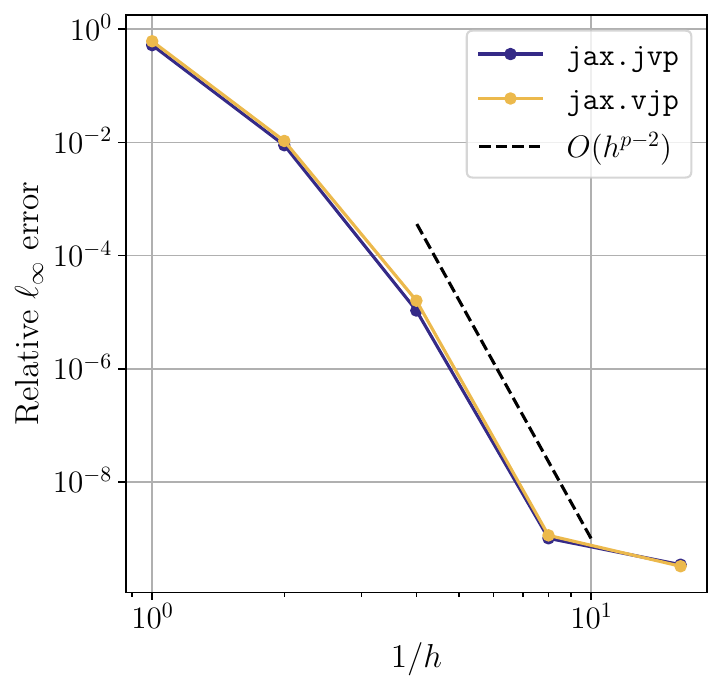}
        \caption{}
        \label{fig:autodiff_accuracy}
    \end{subfigure}%
    \begin{subfigure}[b]{0.3\textwidth}
        \centering
        \includegraphics[height=4cm]{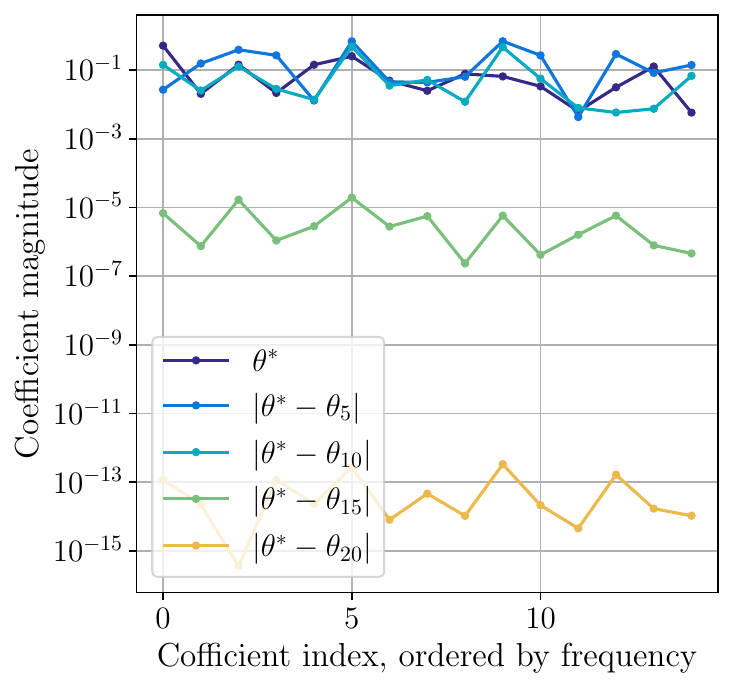}
        \caption{}
        \label{fig:coeffs}
    \end{subfigure}
    \begin{subfigure}[b]{0.3\textwidth}
        \centering
        \includegraphics[height=4cm]{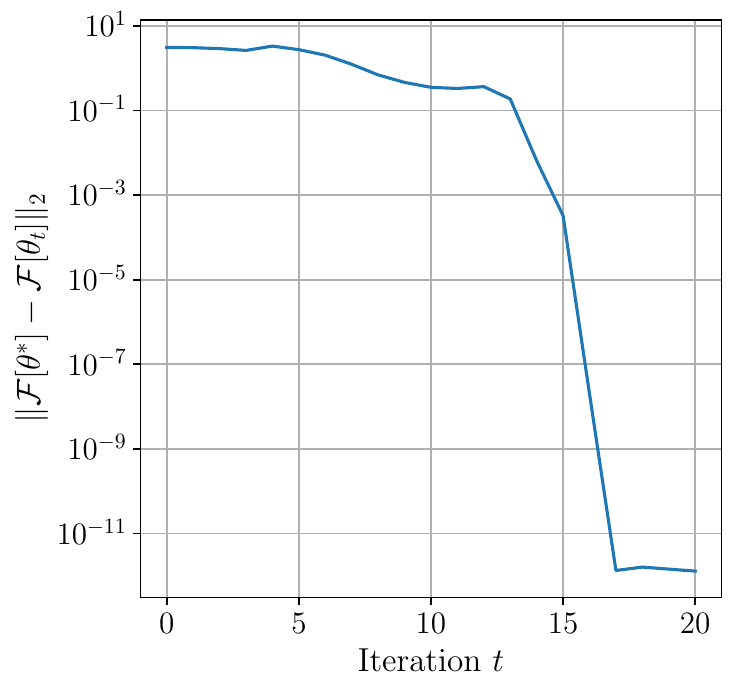}
        \caption{}
        \label{fig:residuals}
    \end{subfigure}
    \caption{
        Our GPU-accelerated PDE solver can interface with automatic differentiation to rapidly solve inverse problems. 
        In \cref{fig:autodiff_accuracy}, \revone{we show that automatic differentiation applied to our HPS method implementing $\cF \circ \cB$ converges at high order to the Fr\'echet derivative of this operator. Experimental details are available in \mbox{\cref{app:adjoints_experimental_details}}. In \mbox{\cref{fig:coeffs}}, we show the magnitude of the ground-truth coefficients $\theta^*$ as well as the component-wise errors of intermediate estimates}. 
        \cref{fig:residuals} shows the objective value of the optimization problem, which reaches \revone{near} machine precision in \revone{17} iterations.}
    \label{fig:inverse_scattering}
\end{figure}

\section{Numerical examples in three dimensions}
\label{sec:3D_examples}
In the three-dimensional examples, we focus on problems with localized regions of high variation. Our adaptive discretization method described in \cref{sec:hardware_acceleration} is designed for such problems. 
In these experiments, we use a single Nvidia H100 GPU with 80GB of on-device memory and 200GB of host RAM.

\subsection{Adaptive refinement on a problem with known solution}
\label{sec:3D_wave_front}
\begin{figure}
    \begin{subfigure}[b]{0.5\textwidth}
        \centering
        \includeinkscape[height=6cm]{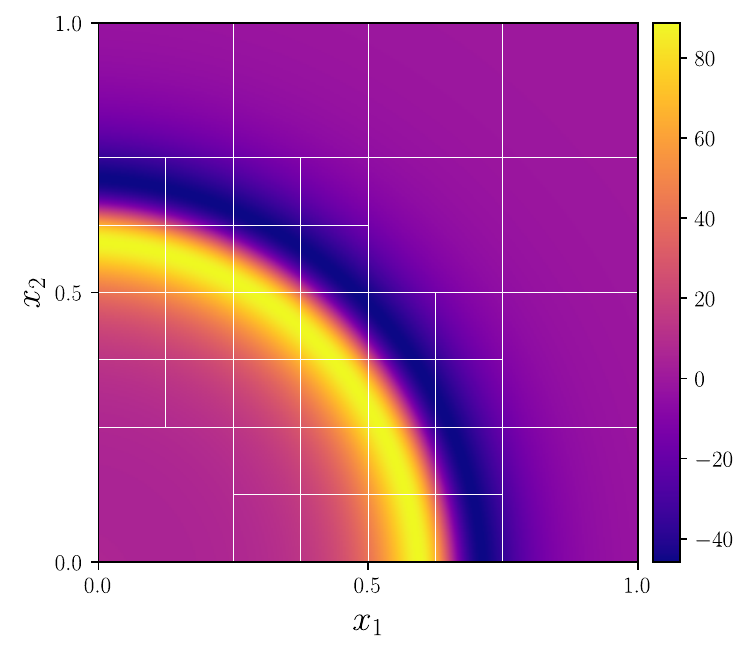}
        \caption{}
        \label{fig:adaptive_source}
    \end{subfigure}
    \begin{subfigure}[b]{0.5\textwidth}
        \centering
        \includeinkscape[height=6cm]{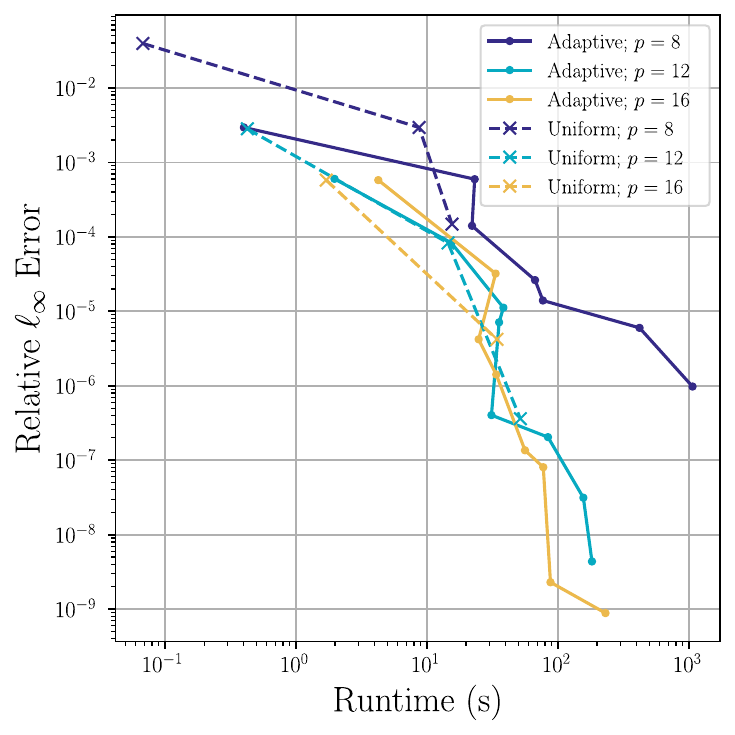}
        \caption{}
        \label{fig:adaptive_vs_uniform_subfig}
    \end{subfigure}
    \caption{Adaptive refinement allows for more accuracy before encountering memory bottlenecks. \cref{fig:adaptive_source} shows the source function restricted to the $x_3=0$ plane and the adaptive mesh formed with tolerance $1\times 10^{-7}$ and Chebyshev parameter $p=16$. In \cref{fig:adaptive_vs_uniform_subfig}, we study the runtime and error of different refinement strategies applied to \cref{problem:3D_wave_front}. For each step of the uniform refinement curve, we refined the uniform grid by one more level. For the adaptive refinement curve, we decreased the adaptive refinement tolerance by a factor of $10$. For all methods, we refined until running out of memory on the GPU during the build stage. }
    \label{fig:adaptive_vs_uniform}
    \vspace{-1em}
\end{figure}

In this example, we study the convergence of our method on a three-dimensional problem with a known analytical solution. 
We build a problem that is solved by a ``wavefront'' located along a three-dimensional curved surface. The problem is given by 
    \begin{align}
    \left\{
        \begin{aligned}
            \Delta u(x) &= f(x), && \quad x \in [0, 1]^3, \\
            u(x) &= g(x), && \quad x \in \partial [0,1]^3.
        \end{aligned}
        \right.
        \label{problem:3D_wave_front}
    \end{align}
We manufacture a source term $f(x)$ so the solution takes the form 
\begin{align}
    u(x) = u(x_1,x_2,x_3) = \arctan \left(10  \sqrt{ (x_1 - 0.5)^2 + (x_2 - 0.5)^2 + (x_3 - 0.5)^2 } - 0.7 \right)
\end{align}
and use samples of this function along the boundary to create our boundary data $g$. 
\cref{fig:adaptive_source} shows that $f$ has a localized region of high variation. 
A uniform discretization strategy cannot adapt to the locality of this problem, but our adaptive discretization strategy can adaptively refine the octree to place a higher density of discretization points in this region. 
The adaptive discretization also uses larger leaves, where possible, on the parts of the domain with a very smooth source function.

In \cref{fig:adaptive_vs_uniform_subfig}, we show accuracy versus runtime for both a uniform and adaptive discretization applied to this problem. For all methods, we increased the number of discretization points until saturating the GPU's memory limit when inverting the highest-level $\boldsymbol{D}$ matrix. 
While the uniform methods are fast on this problem, they are not highly accurate because they cannot use more than $L=3$ levels of uniform refinement. 
Our adaptive method computes solutions with much higher accuracy before saturating the GPU's memory limit by adaptively placing the discretization points in regions where the source and solution have high variation. 
\cref{tab:merge_sizes} shows the size of the highest-level merge matrix $\boldsymbol{D}$ for selected points on this graph.

\subsection{Linearized Poisson--Boltzmann equation}
An example application where the data and solution have local regions of high variation is the linearized Poisson--Boltzmann equation, a model of the electrostatic properties of a molecule in a solution.
This can, for example, be used to compute the stability of a given molecular configuration in a solution.
A standard model, developed in \cite{grant_smooth_2001}, starts with atoms represented by point charges $\{z^{(i)} \}_{i=1}^{N_z}, z^{(i)} \in \R^3$. These atoms give rise to a charge distribution $\rho(x)$,
\begin{align}
    \rho(x) = \sum_{i=1}^{N_z} e^{-\delta \| x - z^{(i)} \|^2 },
\end{align}
and a spatially-varying permittivity function $\varepsilon(x)$,
\begin{align}
    \varepsilon(x) = \epsilon_0 + (\epsilon_\infty - \epsilon_0) e^{-A \rho(x)}.
\end{align}
We use parameters $N_z = 50$, $\delta = 45$, $\epsilon_0 = 16$, $\epsilon_\infty = 100$, and $A = 10$ \citep{grant_smooth_2001}.
The atomic centers $\{z^{(i)} \}_{i=1}^{N_z}$ are drawn uniformly from the box $[-0.5, 0.5]^3$.
We can now model the electrostatic potential $u(x)$ which is implicitly defined by the linearized Poisson--Boltzmann equation:
\begin{align}
\left\{
    \begin{aligned}
        \nabla \cdot \left( \varepsilon(x) \cdot \nabla u(x) \right) &= - \rho(x), && \quad x \in [-1,1]^3, \\
         u(x) &= 0, && \quad x \in \partial [-1,1]^3.
    \end{aligned}
    \right.
    \label{eq:poisson_boltzmann}
\end{align}

Existing approaches, such as a fast integral equation method \citep{vico_fast_2016} and finite difference schemes \citep{nicholls_rapid_1991,colmenares_gpu_2014} solve a simplified version of \cref{eq:poisson_boltzmann}, where the permittivity function is replaced by one derived from van der Waals surfaces:
\newcommand{\vdW}{\text{vdW}}
\begin{align}
    \varepsilon_{\vdW}(x) &= q(x) \epsilon_0 + (1 - q(x)) (\epsilon_{\infty} - \epsilon_0), \\
    q(x) &= 1 - \prod_{i=1}^{N_z}\left[ 1 - e^{-\delta \|x - z^{(i)} \|^2} \right].
\end{align}
$\varepsilon_{\vdW}(x)$ is an easier function to resolve to high accuracy as it lacks the steep gradients observed in $\varepsilon(x)$. However, \cite{grant_smooth_2001} reports ``experimentation with a number of dielectric mapping functions using [$\varepsilon_{\vdW}$] produced dielectric functions that increase toward solvent values far too rapidly with distance from atomic centers,'' and ``[$\varepsilon_{\vdW}$] also produced undesired patches of high dielectric inside proteins.'' We use our GPU-accelerated adaptive HPS method to solve \cref{eq:poisson_boltzmann} with both permittivity models.

\begin{figure}[!ht]
\centering
    \begin{subfigure}[b]{0.3\textwidth}
        \centering
        \includeinkscape[height=4cm]{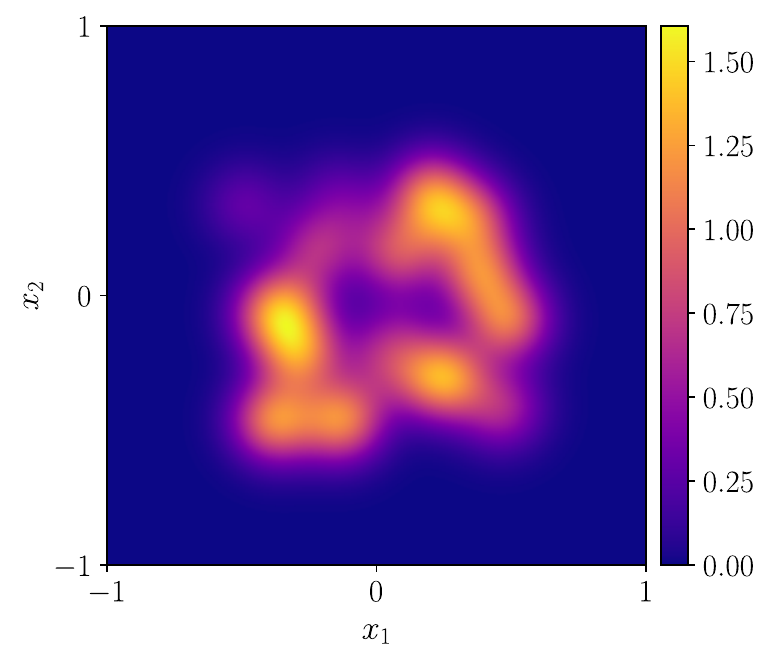}
        \caption{Source $\rho(x)$}
        \label{fig:rho_z}
    \end{subfigure}%
    \begin{subfigure}[b]{0.3\textwidth}
        \centering
        \includeinkscape[height=4cm]{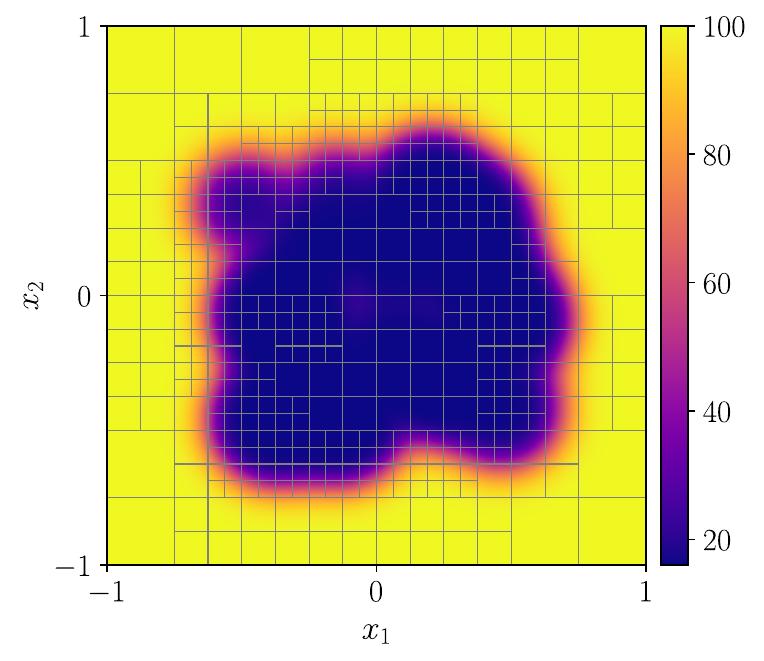}
        \caption{Permittivity $\varepsilon(x)$}
    \end{subfigure}
    \begin{subfigure}[b]{0.3\textwidth}
        \includeinkscape[height=4cm]{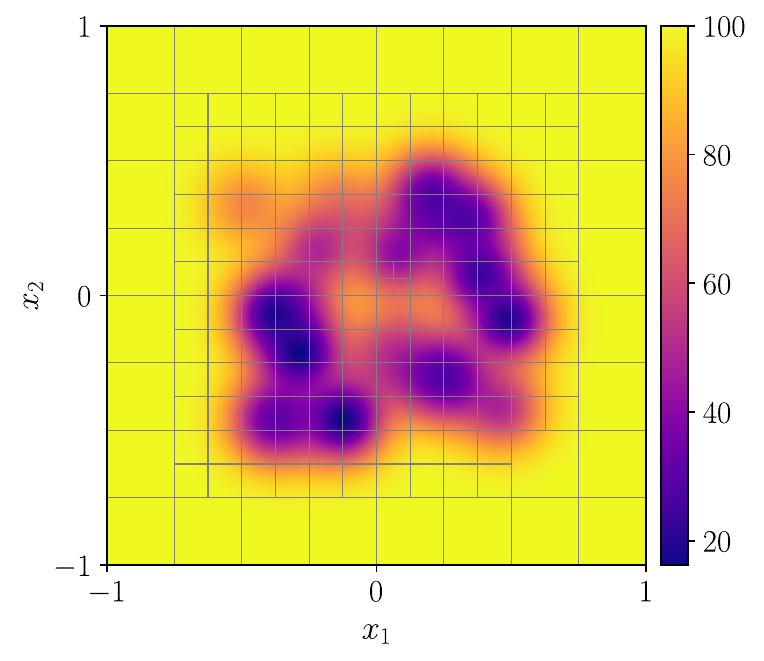}
        \caption{Permittivity $\varepsilon_{\vdW}(x)$}
    \end{subfigure}
    \caption{Visualizing the variable coefficients and source term of our problem. These plots show the source function $\rho(x)$ and permittivity functions $\varepsilon(x)$ and $\varepsilon_{\vdW}(x)$ restricted to the plane $x_3=0$. The adaptive discretizations formed using $p=8$ and tolerance $1\times 10^{-4}$ are shown. This adaptive discretization is found by forming the union of meshes adaptively refined on the source, the permittivity, and components of the gradient of the permittivity. }
    \label{fig:poisson_boltzmann_data}
\end{figure}

To form an adaptive discretization for this problem, we refine a discretization tree given the charge distribution $\rho$, the permittivity $\varepsilon$, and the components of $\nabla \varepsilon$. 
\cref{fig:poisson_boltzmann_data} shows a 2D slice of the charge distribution and permittivity, along with the discretization found by this refinement process. 
\cref{tab:poisson-Boltzmann_results} gives statistics about the discretization and runtime for a range of tolerances and Chebyshev parameters $p$. \revtwo{In this table, $n_{\text{leaves}}$ is the number of leaves of the resulting discretization tree, $N=n_{\text{leaves}}p^3$ is the total number of discretization points, Max Depth is the maximum number of levels of refinement in the discretization tree, and Runtime measures the wall-clock time in seconds to compute the adaptive discretization and execute all parts of the HPS algorithm.}

\begin{table}
    \vspace{-0.5em}
    \centering
    \begin{minipage}{0.47\linewidth}
        \resizebox{\linewidth}{!}{\begin{tabular}{rrrrrr}
            \toprule
            $p$ & \textbf{Tolerance} & $\nleaves$ & $N$ & \textbf{Max Depth} & \textbf{Runtime (s)} \\
            \midrule
            $8$ & $ 10^{-1}$ & $358$ & $183{,}296$ & $3$ & $48.4$ \\
            $8$ & $ 10^{-2}$ & $1{,}436$ & $735{,}232$ & $4$ & $177.8$ \\
            $8$ & $ 10^{-3}$ & $1{,}884$ & $964{,}608$ & $5$ & $218.9$ \\
            $8$ & $ 10^{-4}$ & $7{,}659$ & $3{,}921{,}408$ & $5$ & $1{,}523.0$ \\
            \midrule
            $10$ & $ 10^{-1}$ & $253$ & $253{,}000$ & $3$ & $52.5$ \\
            $10$ & $ 10^{-2}$ & $449$ & $449{,}000$ & $4$ & $76.3$ \\
            $10$ & $ 10^{-3}$ & $1{,}625$ & $1{,}625{,}000$ & $4$ & $243.4$ \\
            $10$ & $ 10^{-4}$ & $1{,}982$ & $1{,}982{,}000$ & $5$ & $319.1$ \\
            \midrule
            $12$ & $ 10^{-1}$ & $99$ & $ 171{,}072$ & $3$ & $41.5$ \\
            $12$ & $ 10^{-2}$ & $386$ & $ 667{,}008$ & $3$ & $91.0$ \\
            $12$ & $ 10^{-3}$ & $715$ & $ 1{,}235{,}520$ & $4$ & $190.6$ \\
            $12$ & $ 10^{-4}$ & $1{,}695$ & $ 2{,}928{,}960$ & $4$ & OOM \\
            \midrule
            $16$ & $ 10^{-1}$ & $64$ & $ 262{,}144$ & $2$ & $53.9$ \\
            $16$ & $ 10^{-2}$ & $204$ & $ 835{,}584$ & $3$ & $136.2$ \\
            $16$ & $ 10^{-3}$ & $365$ & $ 1{,}495{,}040$ & $3$ & OOM \\
            $16$ & $ 10^{-4}$ & $470$ & $ 1{,}925{,}120$ & $4$ & OOM \\
            \bottomrule
        \end{tabular}}
    \end{minipage}
    \begin{minipage}{0.48\linewidth}
        \resizebox{\linewidth}{!}{
        \begin{tabular}{rrrrrr}
            \toprule
            $p$ & \textbf{Tolerance} & $\nleaves$ & $N$ & \textbf{Max Depth} & \textbf{Runtime (s)} \\
            \midrule
            $8$ & $ 1 \times 10^{-3}$ & $1{,}093$ & $ 559{,}616$ & $4$ & $154.4$ \\
            $8$ & $ 1 \times 10^{-4}$ & $1{,}737$ & $ 889{,}344$ & $5$ & $227.2$ \\
            $8$ & $ 1 \times 10^{-5}$ & $5{,}111$ & $ 2{,}616{,}832$ & $5$ & $869.0$ \\
            $8$ & $ 1 \times 10^{-6}$ & $9{,}423$ & $ 4{,}824{,}576$ & $5$ & OOM \\
            \midrule
            $10$ & $ 1 \times 10^{-3}$ & $435$ & $ 435{,}000$ & $4$ & $94.5$ \\
            $10$ & $ 1 \times 10^{-4}$ & $1{,}016$ & $ 1{,}016{,}000$ & $4$ & $193.3$ \\
            $10$ & $ 1 \times 10^{-5}$ & $1{,}485$ & $ 1{,}485{,}000$ & $4$ & $241.7$ \\
            $10$ & $ 1 \times 10^{-6}$ & $2{,}605$ & $ 2{,}605{,}000$ & $5$ & $483.0$ \\
            \midrule
            $12$ & $ 1 \times 10^{-3}$ & $274$ & $ 473{,}472$ & $3$ & $93.7$ \\
            $12$ & $ 1 \times 10^{-4}$ & $477$ & $ 824{,}256$ & $4$ & $145.8$ \\
            $12$ & $ 1 \times 10^{-5}$ & $974$ & $ 1{,}683{,}072$ & $4$ & $275.2$ \\
            $12$ & $ 1 \times 10^{-6}$ & $1{,}366$ & $ 2{,}360{,}448$ & $4$ & OOM \\
            \midrule
            $16$ & $ 1 \times 10^{-3}$ & $127$ & $ 520{,}192$ & $3$ & $119.7$ \\
            $16$ & $ 1 \times 10^{-4}$ & $239$ & $ 978{,}944$ & $3$ & $172.2$ \\
            $16$ & $ 1 \times 10^{-5}$ & $323$ & $ 1{,}323{,}008$ & $3$ & $227.0$ \\
            $16$ & $ 1 \times 10^{-6}$ & $456$ & $ 1{,}867{,}776$ & $4$ & OOM \\
            \bottomrule
        \end{tabular}
        }
    \end{minipage} 
\caption{Resource usage statistics for using our 3D adaptive HPS method applied to the linearized Poisson--Boltzmann equation (\cref{eq:poisson_boltzmann}) using permittivity $\varepsilon(x)$ (left) and simplified permittivity $\varepsilon_{\vdW}(x)$ (right).
We use ``OOM'' to indicate which discretizations caused out of memory errors when invering the final $\boldsymbol{D}$ matrix.}
\label{tab:poisson-Boltzmann_results}
\vspace{-1em}
\end{table}

\section{Conclusion}
\label{sec:conclusion}
This paper presents methods for efficiently accelerating HPS algorithms using general-purpose GPUs. 
Because there is a large amount of inherent parallelism in the structure of the HPS algorithms, they are a natural target for GPU acceleration once adjustments are made to reduce memory complexity. 
We introduce methods for reducing the memory footprint of HPS algorithms for problems in two and three dimensions.

This work leaves open important questions and avenues for improvement.
While our method can efficiently interface with automatic differentiation, we could, in principle, gain much more efficiency by implementing custom automatic differentiation rules to reuse the precomputed solution operators, like those derived in \citet{borges_high_2016}. 
Our methods of reducing memory complexity could be pushed further by using a hybrid approach, i.e., by performing a few levels of merging via dense linear algebra and then relying on a sparse direct solver such as \citet{yesypenko_slablu_2024} or \citet{Kump_Yesypenko_Martinsson_2025} for higher-level merges. 
We hypothesize such a hybrid approach would greatly reduce runtimes for very large problems by reducing the ranks needed to accurately resolve the sparse system matrix. 
Finally, our methods could be extended to unstructured meshes and surfaces, a setting where nonuniform merge operations make parallel implementations challenging.

\section*{Author Contributions}
\begin{itemize}
\item 
\textbf{Owen Melia:} Conceptualization; Methodology; Software; Writing - original draft. 
\item \textbf{Daniel Fortunato:} Conceptualization; Methodology; Supervision; Writing - review \& editing. 
\item \textbf{Jeremy Hoskins:} Conceptualization; Methodology; Supervision; Writing - review \& editing. 
\item \textbf{Rebecca Willett:} Conceptualization; Methodology; Supervision; Funding acquisition; Project administration; Writing - review \& editing. 
\end{itemize}

\section*{Data Availability}
Our open-source JAX implementation is publicly available at \url{https://github.com/meliao/jaxhps}.

\section*{Acknowledgment}
The authors would like to thank Manas Rachh, Leslie Greengard, and Olivia Tsang for many useful discussions. 
The authors are grateful to the Flatiron Institute for providing the computational resources used to conduct the experiments in this work. 
OM and RW gratefully acknowledge the support of NSF DMS-2023109, DOE DE-SC0022232, the Physics Frontier Center for Living Systems funded by the National Science Foundation (PHY-2317138), and the support of the Margot and Tom Pritzker Foundation.
The Flatiron Institute is a division of the Simons Foundation.

\bibliographystyle{elsarticle-harv} 
\bibliography{ref}

\appendix
\setcounter{MaxMatrixCols}{20}
\section{Full algorithms for 2D problems using DtN matrices}
\label{appendix:2D_DtN}
In this section, we describe the details for the two-dimensional version of our method which merges DtN matrices. 
In this version of the algorithm, the outgoing boundary data $\boldsymbol{h}$ tabulates the outward-pointing normal derivative of the particular solution, and the incoming boundary data $\boldsymbol{g}$ tabulates the homogenous solution values restricted to patch boundaries. $\boldsymbol{T}$ is a Dirichlet-to-Neumann matrix. 

In this section, we use $\boldsymbol{I}_{a\times a}$ to denote the identity matrix of shape $a \times a$ and $\boldsymbol{0}_{d}$ to denote a length-$d$ vector filled with $0$'s. 
When defining matrices blockwise, we use $\boldsymbol{0}$ to denote a block filled with $0$'s, and assume the shape of the block can be determined from the nonzero blocks sharing the same rows and columns.

\subsection{Local solve stage}
Recall from \cref{sec:discretization} that we use a tensor product of order-$p$ Chebyshev--Lobatto points to discretize the interior of each leaf. This results in a grid with $p^2$ discretization points, and $4p-4$ of these points lie on the boundary of the leaf. 
We use order-$q$ Gauss--Legendre points to discretize each side of the leaf's boundary, so there are $4q$ boundary points in total. 
Thus, to translate the information between the interior and boundary of each leaf, we need to compute spectral differentiation matrices and matrices interpolating between the $p$ Chebyshev and $q$ Gauss points. 
In particular, we need to precompute the following matrices:
\begin{itemize}
    \item $\boldsymbol{P}$, with shape $4p-4 \times 4q$, is the operator mapping data sampled on the Gauss boundary points to data sampled on the $4p-4$ Chebyshev points located on the boundary of the leaf.
    This matrix is constructed using a barycentric Lagrange interpolation matrix mapping from Gauss to Chebyshev points on one side of the leaf; this interpolation matrix is repeated for the other sides. Rows corresponding to the Chebyshev points on the corners of the leaf average the contribution from the two adjoining panels. 
    \item $\boldsymbol{Q}$, with shape  $4q \times p^2$, performs spectral differentiation on the $p^2$ Chebyshev points followed by interpolation to the Gauss boundary points.
    This matrix is formed by stacking the relevant rows of Chebyshev spectral differentiation matrices to form an operator which evaluates normal derivatives on the $4p-4$ boundary Chebyshev points, and then composing this differentiation operator with a matrix formed from barycentric Lagrange interpolation matrix blocks. 
    These interpolation matrices each map from one Chebyshev panel to one Gauss panel.
\end{itemize}
To work with $\boldsymbol{L}^{(i)}$, the discretization of the differential operator on leaf $i$, it is useful to identify $I_i$ and $I_e$, the sets of discretization points corresponding to the $(p-2)^2$ interior and  $4p-4$ exterior Chebyshev points, respectively. 
Now, we can fully describe the local solve stage in \cref{alg:2D_DtN_local}.

\begin{algorithm}[h!]
    \caption{
       2D DtN local solve stage.
    }
    \label{alg:2D_DtN_local}
    \KwIn{Discretized differential operators $\{\boldsymbol{L}^{(i)} \}_{i=1}^{\nleaves}$; 
    discretized source vectors $\{ \boldsymbol{f}^{(i)} \}_{i=1}^{\nleaves}$; 
    precomputed interpolation and differentiation matrices $\boldsymbol{P}$ and $\boldsymbol{Q}$
    }
    \For(){$i=1, \hdots , \nleaves$} {
        Invert $\boldsymbol{L}^{(i)}(I_i , I_i)$ \\
        $\boldsymbol{Y}^{(i)} = \begin{bmatrix}
            \boldsymbol{I}_{4p-4 \times 4p-4} \\
            - \left( \boldsymbol{L}^{(i)}(I_i , I_i)\right)^{-1} \boldsymbol{L}^{(i)}(I_i , I_e)
        \end{bmatrix} \boldsymbol{P}$ \label{eqn:def_Y_DtN} \\
        $\boldsymbol{v}^{(i)} = \begin{bmatrix}
            \boldsymbol{0}_{4p-4} \\
            - \left( \boldsymbol{L}^{(i)}(I_i , I_i)\right)^{-1} \boldsymbol{f}^{(i)}\left( I_i \right)
        \end{bmatrix}$ \label{eqn:dev_v_DtN} \\
        $\boldsymbol{T}^{(i)} = \boldsymbol{Q} \boldsymbol{Y}^{(i)} $ \\
        $\boldsymbol{h}^{(i)} = \boldsymbol{Q} \boldsymbol{v}^{(i)}$
    }
    \KwResult{Poincar\'e--Steklov matrices $\{ \boldsymbol{T}^{(i)} \}_{i=1}^{\nleaves}$; 
    outgoing boundary data $\{ \boldsymbol{h}^{(i)} \}_{i=1}^{\nleaves}$; 
    interior solution matrices $\{ \boldsymbol{Y}^{(i)} \}_{i=1}^{\nleaves}$, 
    leaf-level particular solutions $\{ \boldsymbol{v}^{(i)} \}_{i=1}^{\nleaves}$ }
\end{algorithm}

\subsection{Merge stage}
\label{appendix:2D_DtN_merge}
\begin{figure}
    \centering
    
\begin{tikzpicture}
    \def\s{2}
    \def\offset{0.1}

    \node at (\s/2, \s/2) {\large $\Omega_a$}; %
    \node at (3*\s/2, \s/2) {\large $\Omega_b$}; %
    \node at (\s/2, 3*\s/2) {\large $\Omega_d$}; %
    \node at (3*\s/2, 3*\s/2) {\large $\Omega_c$}; %

    \draw[thick] (-\offset, -\offset) -- (\s - \offset, -\offset);
    \draw[thick] (-\offset, -\offset) -- (-\offset, \s - \offset);
    \node[anchor=north east] at (-\offset, -\offset) {\large 1};

    \draw[thick] (\s + \offset, -\offset) -- (2*\s + 2*\offset, -\offset);
    \draw[thick] (2*\s + 2*\offset, -\offset) -- (2*\s + 2*\offset, \s - \offset);
    \node[anchor=north west] at (2*\s + 2*\offset, -\offset) {\large 2};

    \draw[thick] (2*\s + 2*\offset, \s + \offset) -- (2*\s + 2*\offset, 2*\s + 2*\offset);
    \draw[thick] (\s + \offset, 2*\s + 2*\offset) -- (2*\s + 2*\offset, 2*\s + 2*\offset);
    \node[anchor=south west] at (2*\s + 2*\offset, 2*\s + 2*\offset) {\large 3};

    \draw[thick] (-\offset, \s + \offset) -- (-\offset, 2*\s + 2*\offset);
    \draw[thick] (-\offset, 2*\s + 2*\offset) -- (\s - \offset, 2*\s + 2*\offset);
    \node[anchor=south east] at (-\offset, 2*\s + 2*\offset) {\large 4};

    \draw[thick] (\s, 0) -- (\s, \s - \offset) node[midway, left] {\large 5};
    \draw[thick] (\s + \offset, \s) -- (2*\s, \s) node[midway, below] {\large 6};
    \draw[thick] (\s, \s + \offset) -- (\s, 2*\s) node[midway, right] {\large 7};
    \draw[thick] (\s - \offset, \s) -- (0, \s) node[midway, above] {\large 8};

\end{tikzpicture}
\caption{Visualizing boundary elements $1$ through $8$ for two-dimensional merges.}
\label{fig:2D_merges}
\end{figure}
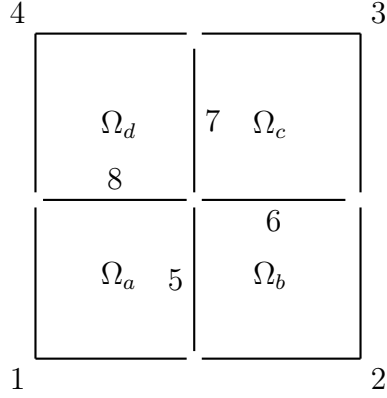
\newcommand{\nside}{n_{\text{side}}}
In the 2D merge stage, we are merging four nodes $\Omega_a, \Omega_b, \Omega_c$, and $\Omega_d$, which have \emph{exterior} and \emph{interior} discretization points. We label the exterior boundary sections $1,2,3$, and $4$, and we label the interior boundary sections $5, 6, 7$, and $8$. See \cref{fig:2D_merges} for a diagram of the different boundary parts. 
Because the merge stage operates completely on data discretized using Gauss--Legendre panels, there are no discretization points at the corners of nodes. This means each discretization point belongs to exactly one part of the boundary. 
During this stage of the algorithm, we will be indexing rows and columns of the Dirichlet-to-Neumann matrices according to these boundary sections. 
For example, we use $\boldsymbol{T}^{(a)}_{1,5}$ to indicate the submatrix of node $a$'s DtN matrix which maps from boundary section $5$ to boundary section $1$.
Suppose each side of $\Omega_a, \Omega_b, \Omega_c$, and $\Omega_d$ is discretized with $\nside$ discretization points; in this case $\boldsymbol{T}^{(a)}_{1,5}$ will have shape $2\nside \times \nside$.

To implement the merge stage, we use sets of constraints to solve for a mapping from given $\boldsymbol{g}_{\ext}$ to unknown $\boldsymbol{g}_{\intt}$. These are vectors tabulating the homogeneous solution along the exterior and interior boundary parts. First are constraints specifying that the solution to the PDE is continuous:
\begin{align}
    \boldsymbol{u}_{\ext} &= \boldsymbol{A} \boldsymbol{g}_{\ext} +  \boldsymbol{B}\boldsymbol{g}_{\intt} + \boldsymbol{h}^{(\child)}_{\ext}.
    \label{eq:2D_DtN_merge_1}
\end{align}
In this equation $\boldsymbol{u}_{\ext}$ is interpreted as the outward-pointing normal derivative of the solution to the PDE restricted to boundary elements $1,2,3,$ and $4$ with boundary data specified by $\boldsymbol{g}_{\ext}$.
In this set of constraints, we define:
\newcommand{\bz}{\boldsymbol{0}}
\begin{align}
    \boldsymbol{h}^{(\child)}_{\ext} &= \begin{bmatrix}
        \boldsymbol{h}^{(a)}_{1} \\
        \boldsymbol{h}^{(b)}_{2} \\
        \boldsymbol{h}^{(c)}_{3} \\
        \boldsymbol{h}^{(d)}_{4} \\
    \end{bmatrix}, \label{eq:2D_DtN_h_child_ext} \\
    \boldsymbol{A} &= \begin{bmatrix}
        \boldsymbol{T}^{(a)}_{1,1} & \bz & \bz & \bz \\ 
        \bz & \boldsymbol{T}^{(b)}_{2,2} & \bz & \bz  \\ 
        \bz & \bz & \boldsymbol{T}^{(c)}_{3,3} &  \bz \\
        \bz & \bz &  \bz & \boldsymbol{T}^{(d)}_{4,4}
    \end{bmatrix}, \label{eq:2D_DtN_A} \\
    \boldsymbol{B} &= \begin{bmatrix}
        \boldsymbol{T}^{(a)}_{1,5} & \bz & \bz & \boldsymbol{T}^{(a)}_{1,8} \\ 
        \boldsymbol{T}^{(b)}_{2,5} & \boldsymbol{T}^{(b)}_{2,6} &\bz & \bz  \\ 
        \bz & \boldsymbol{T}^{(c)}_{3,6} & \boldsymbol{T}^{(c)}_{3,7} &\bz \\ 
        \bz  & \bz & \boldsymbol{T}^{(d)}_{4,7} & \boldsymbol{T}^{(d)}_{4,8}
    \end{bmatrix}. \label{eq:2D_DtN_B}
\end{align}
A second set of constraints enforces that the outward-pointing normal derivatives from neighboring nodes should sum to zero, which is equivalent to enforcing the continuity of the first derivative along the merge interfaces in their respective Cartesian directions. To that end, we use constraints $\boldsymbol{u}^{(a)}_{5} + \boldsymbol{u}^{(b)}_{5} = \boldsymbol{0}_{\nside}$, $\boldsymbol{u}^{(b)}_{6} + \boldsymbol{u}^{(c)}_{6} = \boldsymbol{0}_{\nside}$, and so on. This gives us an equation:
\begin{align}
    \boldsymbol{0}_{4\nside} &= \boldsymbol{C} \boldsymbol{g}_{\ext} +  \boldsymbol{D}\boldsymbol{g}_{\intt} + \boldsymbol{h}^{(\child)}_{\intt}.
    \label{eq:2D_DtN_merge_2}
\end{align} 
In \cref{eq:2D_DtN_merge_2}, we define
\begin{align}
    \boldsymbol{h}_{\intt}^{(\child)} &= \begin{bmatrix}
        \boldsymbol{h}^{(a)}_{5} + \boldsymbol{h}^{(b)}_{5} \\
        \boldsymbol{h}^{(b)}_{6} + \boldsymbol{h}^{(c)}_{6} \\
        \boldsymbol{h}^{(c)}_{7} + \boldsymbol{h}^{(d)}_{7} \\
        \boldsymbol{h}^{(d)}_{8} + \boldsymbol{h}^{(a)}_{8} \\
    \end{bmatrix}, \label{eq:2D_DtN_h_child_int} \\
    \boldsymbol{C} &= \begin{bmatrix}
        \boldsymbol{T}^{(a)}_{5,1} & \boldsymbol{T}^{(b)}_{5,2} & \bz & \bz \\
        \bz & \boldsymbol{T}^{(b)}_{6,2} &  \boldsymbol{T}^{(c)}_{6,3} & \bz \\
        \bz & \bz &  \boldsymbol{T}^{(c)}_{7,3} & \boldsymbol{T}^{(d)}_{7,4} \\
        \boldsymbol{T}^{(a)}_{8,1} & \bz & \bz & \boldsymbol{T}^{(d)}_{8,4} \\
    \end{bmatrix}, \label{eq:2D_DtN_C} \\
    \boldsymbol{D} &= \begin{bmatrix}
        \boldsymbol{T}^{(a)}_{5,5} + \boldsymbol{T}^{(b)}_{5,5} & \boldsymbol{T}^{(b)}_{5,6} & \bz & \boldsymbol{T}^{(a)}_{5,8} \\
        \boldsymbol{T}^{(b)}_{6,5} & \boldsymbol{T}^{(b)}_{6,6} + \boldsymbol{T}^{(c)}_{6,6} & \boldsymbol{T}^{(c)}_{6,7} & \bz \\
        \bz &  \boldsymbol{T}^{(c)}_{7,6} &  \boldsymbol{T}^{(c)}_{7,7} +  \boldsymbol{T}^{(d)}_{7,7} &  \boldsymbol{T}^{(d)}_{7,8} \\
        \boldsymbol{T}^{(a)}_{8,5} & \bz & \boldsymbol{T}^{(d)}_{8,7} & \boldsymbol{T}^{(d)}_{8,8} + \boldsymbol{T}^{(a)}_{8,8}
    \end{bmatrix}. \label{eq:2D_DtN_D}
\end{align}
Now that the matrices and vectors are defined, we can construct the linear system in \cref{eq:merge_lin_system} and compute the merged data.

\section{Full algorithms for 2D problems using ItI matrices}
\label{appendix:2D_ItI}
In this section, we describe the details for the two-dimensional version of our method which merges ItI matrices. 
In this version of the algorithm, the outgoing boundary data $\boldsymbol{h}$ tabulates the outgoing impedance data due to the particular solution, and the incoming boundary data $\boldsymbol{g}$ tabulates incoming impedance data due to the homogeneous solution. $\boldsymbol{T}$ is an impedance-to-impedance matrix. To define the impedance data, we need to choose a value $\eta \in \R_{+}$. In the wave scattering context, we often choose $\eta = k$ \citep{gillman_spectrally_2015}.

As in the previous section, we use  $\boldsymbol{I}_{a\times a}$ to denote the identity matrix of shape $a \times a$ and $\boldsymbol{0}_{d}$ to denote a length-$d$ vector filled with $0$'s. 
We also use $\boldsymbol{0}$ to denote a block filled with $0$'s, and assume the shape of this block can be inferred from the context.
\subsection{Local solve stage}
As before, there are $4p-4$ Chebyshev points on the boundary and $4q$ Gauss points on the boundary, and we must map between these two sets of discretization points. 
In particular, we need to precompute the following matrices:
\begin{itemize}
    \item $\boldsymbol{P}$, with shape $4p-4 \times 4q$, is the operator mapping data sampled on the Gauss boundary points to data sampled on the $4p-4$ Chebyshev points located on the boundary of the leaf.
    This matrix is constructed using a barycentric Lagrange interpolation matrix mapping from Gauss to Chebyshev points on one side of the leaf with the final row deleted; this interpolation matrix is repeated for the other sides.
    \item $\boldsymbol{Q}$, with shape $4q \times 4p$, is the operator mapping data sampled on the Chebyshev points located on the boundary of the leaf to the Gauss points on the boundary of the leaf. This matrix is block-diagonal with four copies of a barycentric Lagrange interpolation matrix mapping from Chebyshev to Gauss points on one side of the leaf. Note this matrix double-counts the Chebyshev points at the corners of the leaf.
    \item $\boldsymbol{N}$, with shape $4p \times p^2$, is an operator mapping from interior solutions to outward-pointing normal derivative data evaluated on the boundary Chebyshev points. Note this operator counts each corner point twice. This matrix is formed by stacking relevant rows of Chebyshev spectral differentiation matrices.
    \item $\boldsymbol{\tilde{N}}$, with shape $4p-4 \times p^2$, is an operator mapping from interior solutions to outward-pointing normal data evaluated on the boundary Chebyshev points. Note this operator counts each corner point only once. This matrix is formed by stacking relevant rows of Chebyshev spectral differentiation matrices.
    \item $\boldsymbol{H}$, with shape $4p \times p^2$, is an operator mapping from interior solutions to evaluations of outgoing impedance data on the Chebyshev boundary discretization points. This matrix is constructed by taking $\boldsymbol{N}$ and subtracting $i\eta \boldsymbol{I}_{p\times p}$ from the appropriate submatrices. Again, this matrix double-counts the Chebyshev discretization points at the corners of the leaf.
    \item $\boldsymbol{G}$, with shape $4p-4 \times p^2$, is an operator mapping from interior solutions to evaluations of incoming impedance data on the Chebyshev boundary discretization points. This matrix is constructed by taking $\boldsymbol{\tilde{N}}$ and adding $i\eta \boldsymbol{I}_{p-1 \times p-1}$ to the appropriate submatrices. 
\end{itemize}
To work with $\boldsymbol{L}^{(i)}$, the discretization of the differential operator on leaf $i$, it is useful to identify $I_i$ and $I_e$, the sets of discretization points corresponding to the $(p-2)^2$ interior and  $4p-4$ exterior Chebyshev points, respectively. As in \cref{appendix:2D_DtN}, we solve the local problem by using these precomputed operators to enforce the differential operator on the interior discretization points and the boundary condition on the boundary discretization points. In this case, the boundary condition is an incoming impedance condition, also known as a Robin boundary condition. Now, we can fully describe the local solve stage in \cref{alg:2D_ItI_local}.

\begin{algorithm}[h!]
    \caption{
       2D ItI local solve stage.
    }
    \label{alg:2D_ItI_local}
    \KwIn{Discretized differential operators $\{\boldsymbol{L}^{(i)} \}_{i=1}^{\nleaves}$; discretized source functions $\{ \boldsymbol{f}^{(i)} \}_{i=1}^{\nleaves}$; precomputed interpolation and differentiation matrices $\boldsymbol{P}, \boldsymbol{Q}, \boldsymbol{H}$, and $\boldsymbol{G}$.
    }
    \For(){$i=1, \hdots , \nleaves$} {
        $\boldsymbol{B}^{(i)} = \begin{bmatrix}
            \boldsymbol{G} \\ \boldsymbol{L}^{(i)}\left( I_i , :\right)
        \end{bmatrix}$ \\
        Invert $\boldsymbol{B}^{(i)}$ \\
        $\boldsymbol{Y}^{(i)} = \left( \boldsymbol{B}^{(i)} \right)^{-1}\left(:, I_e\right) \boldsymbol{P}$ \\
        $\boldsymbol{v}^{(i)} = \left( \boldsymbol{B}^{(i)} \right)^{-1}\left( :, I_i \right) \boldsymbol{f}^{(i)}\left(I_i\right)$ \\
        $\boldsymbol{T}^{(i)} = \boldsymbol{Q} \boldsymbol{H} \boldsymbol{Y}^{(i)} $ \\
        $\boldsymbol{h}^{(i)} = \boldsymbol{Q} \boldsymbol{H} \boldsymbol{v}^{(i)}$
    }
    \KwResult{Poincar\'e--Steklov matrices $\{ \boldsymbol{T}^{(i)} \}_{i=1}^{\nleaves}$; outgoing boundary data $\{ \boldsymbol{h}^{(i)} \}_{i=1}^{\nleaves}$; interior solution matrices $\{ \boldsymbol{Y}^{(i)} \}_{i=1}^{\nleaves}$, leaf-level particular solutions $\{ \boldsymbol{v}^{(i)} \}_{i=1}^{\nleaves}$ }
\end{algorithm}

\subsection{Merge stage}
\label{appendix:2D_ItI_merge}

As in \cref{appendix:2D_DtN}, we are merging four nodes $\Omega_a, \Omega_b, \Omega_c$, and $\Omega_d$ with boundary parts labeled $1,2,\dots, 8$. See \cref{fig:2D_merges} for a diagram of the different boundary parts. 

To implement the merge stage, we use sets of constraints to solve for a mapping from given $\boldsymbol{g}_{\ext}$ to unknown $\boldsymbol{g}_{\intt}$. These are vectors tabulating incoming impedance data due to the homogeneous solution along the boundary parts. 
Because $\boldsymbol{g}_{\intt}$ tabulates impedance data, we must represent the incoming data with respect to neighboring nodes separately. For example, we must represent $\boldsymbol{g}_{5}^{(a)}$, the data along boundary element $5$ incoming to node $a$, separately from $\boldsymbol{g}_{5}^{(b)}$. To that end, we use $\boldsymbol{g}_{\intt} = \left[ \boldsymbol{g}_{5}^{(a)},  \boldsymbol{g}_{8}^{(a)}, \boldsymbol{g}_{6}^{(c)}, \boldsymbol{g}_{7}^{(c)}, \boldsymbol{g}_{5}^{(b)}, \boldsymbol{g}_{6}^{(b)}, \boldsymbol{g}_{7}^{(d)}, \boldsymbol{g}_{8}^{(d)} \right]^\top$. The ordering of these boundary elements is chosen specifically to reduce the computation during the merge, which will become apparent later. Again, we use $\nside$ to denote the number of discretization points along each side of the nodes being merged, so $\boldsymbol{g}_{\intt}$ has length $8 \nside$.

The first set of constraints specify that the solution to the PDE is continuous:
\begin{align}
    \boldsymbol{u}_{\ext} &= \boldsymbol{A} \boldsymbol{g}_{\ext} +  \boldsymbol{B}\boldsymbol{g}_{\intt} + \boldsymbol{h}^{(\child)}_{\ext}.
    \label{eq:2D_ItI_merge_1}
\end{align}
In this equation $\boldsymbol{u}_{\ext}$ is interpreted as the outgoing impedance data of the solution to the PDE restricted to the merged nodes with boundary data specified by $\boldsymbol{g}_{\ext}$.
In this set of constraints, we define:

\begingroup
\allowdisplaybreaks
\begin{align}
    \boldsymbol{h}^{(\child)}_{\ext} &= \begin{bmatrix}
        \boldsymbol{h}^{(a)}_{1} \\
        \boldsymbol{h}^{(b)}_{2} \\
        \boldsymbol{h}^{(c)}_{3} \\
        \boldsymbol{h}^{(d)}_{4} \\
    \end{bmatrix}, \\
    \boldsymbol{A} &= \begin{bmatrix}
        \boldsymbol{T}^{(a)}_{1,1} & \bz & \bz & \bz \\ 
        \bz & \boldsymbol{T}^{(b)}_{2,2} & \bz & \bz  \\ 
        \bz & \bz & \boldsymbol{T}^{(c)}_{3,3} &  \bz \\
        \bz & \bz &  \bz & \boldsymbol{T}^{(d)}_{4,4}
    \end{bmatrix}, \\
    \boldsymbol{B} &= \begin{bmatrix}
        \boldsymbol{T}^{(a)}_{1,5} & \boldsymbol{T}^{(a)}_{1,8}  & \bz & \bz & \bz &  \bz & \bz & \bz \\ 
        \bz & \bz & \bz & \bz & \boldsymbol{T}^{(b)}_{2,5} & \boldsymbol{T}^{(b)}_{2,6} &\bz & \bz \\ 
        \bz & \bz & \boldsymbol{T}^{(c)}_{3,6} & \boldsymbol{T}^{(c)}_{3,7} &\bz & \bz & \bz & \bz   \\ 
        \bz & \bz & \bz & \bz & \bz & \bz & \boldsymbol{T}^{(d)}_{4,7} & \boldsymbol{T}^{(d)}_{4,8} 
    \end{bmatrix}.
\end{align}
\endgroup
A second set of constraints specifies that the outgoing total solution's impedance data from one node must be opposite to the incoming homogeneous solution's impedance data for the neighboring node.
For example, along merge interface $5$, we enforce this constraint:
\begin{align}
     \bz_{\nside} &= \boldsymbol{u}^{(b)}_5 + \boldsymbol{g}^{(a)}_5.
\end{align}
In this equation, $\boldsymbol{u}^{(b)}_5$ is the outgoing impedance data due to the total solution of the PDE restricted to the merged nodes with boundary condition $\boldsymbol{g}_{\ext}$. 
The normal derivative is oriented relative to node $b$. We can expand $\boldsymbol{u}^{(b)}_5$ to find:
\begin{align}
    \bz_{\nside} &= \boldsymbol{g}^{(a)}_5 + \boldsymbol{T}_{5,5}^{(a)} \boldsymbol{g}^{(a)}_5 + \boldsymbol{T}_{5,8}^{(a)}\boldsymbol{g}^{(a)}_8 + \boldsymbol{T}_{5,5}^{(b)} \boldsymbol{g}^{(b)}_5 + \boldsymbol{T}_{5,6}^{(b)}\boldsymbol{g}^{(b)}_6 + \boldsymbol{T}_{5,5}^{(a)} \boldsymbol{g}^{(a)}_5 + \boldsymbol{T}_{5,1}^{(a)}\boldsymbol{g}^{(a)}_1 + \boldsymbol{h}^{(a)}_5.
\end{align}

Similar equalities hold in each direction along each merge interface. We can expand to form a second system of constraints:
\begin{align}
    -\boldsymbol{h}_{\intt}^{(\child)} &= \boldsymbol{C} \boldsymbol{g}_{\ext} +  \boldsymbol{D}\boldsymbol{g}_{\intt},
    \label{eq:2D_ItI_merge_2}
\end{align} 
where we define
\begingroup
\allowdisplaybreaks
\begin{align}
    \boldsymbol{h}_{\intt}^{(\child)} &= \begin{bmatrix}
        \boldsymbol{h}^{(b)}_{5} \\
        \boldsymbol{h}^{(d)}_{8} \\
        \boldsymbol{h}^{(b)}_{6} \\
        \boldsymbol{h}^{(d)}_{7} \\
        \boldsymbol{h}^{(a)}_{5} \\
        \boldsymbol{h}^{(c)}_{6} \\
        \boldsymbol{h}^{(c)}_{7} \\
        \boldsymbol{h}^{(a)}_{8} \\
    \end{bmatrix}, \\
    \boldsymbol{C} &= \begin{bmatrix}
        \bz & \boldsymbol{T}^{(b)}_{5,2} &  \bz & \bz \\
        \bz & \bz & \bz & \boldsymbol{T}^{(d)}_{8,4} \\
        \bz &  \boldsymbol{T}^{(b)}_{6,2} & \bz & \bz \\ 
        \bz & \bz & \bz & \boldsymbol{T}^{(d)}_{7,4} \\
        \boldsymbol{T}^{(a)}_{5,1} & \bz & \bz & \bz \\
        \bz & \bz &  \boldsymbol{T}^{(c)}_{6,3} & \bz \\
        \bz & \bz & \boldsymbol{T}^{(c)}_{7,3} & \bz \\ 
        \boldsymbol{T}^{(a)}_{8,1} & \bz & \bz &  \bz \\ 
    \end{bmatrix}, \\
    \boldsymbol{D} &= \boldsymbol{I}_{8 \nside \times 8 \nside} + \begin{bmatrix}
        \bz  & \bz & \bz & \bz & \boldsymbol{T}^{(b)}_{5,5}  &  \boldsymbol{T}^{(b)}_{5,6} & \bz & \bz \\
        \bz & \bz & \bz & \bz & \bz & \bz & \boldsymbol{T}^{(d)}_{8,7} & \boldsymbol{T}^{(d)}_{8,8}  \\
        \bz  & \bz & \bz & \bz & \boldsymbol{T}^{(b)}_{6,5}  &  \boldsymbol{T}^{(b)}_{6,6} & \bz & \bz \\
        \bz & \bz & \bz & \bz & \bz & \bz & \boldsymbol{T}^{(d)}_{7,7} & \boldsymbol{T}^{(d)}_{7,8}   \\
        \boldsymbol{T}^{(a)}_{5,5} & \boldsymbol{T}^{(a)}_{5,8} &  \bz & \bz & \bz & \bz & \bz & \bz  \\
        \bz & \bz & \boldsymbol{T}^{(c)}_{6,6} & \boldsymbol{T}^{(c)}_{6,7} & \bz & \bz & \bz & \bz \\
        \bz & \bz & \boldsymbol{T}^{(c)}_{7,6} & \boldsymbol{T}^{(c)}_{7,7} & \bz & \bz & \bz & \bz \\
        \boldsymbol{T}^{(a)}_{8,5} & \boldsymbol{T}^{(a)}_{8,8} &  \bz & \bz & \bz & \bz & \bz & \bz \label{eq:ItI_D}  \\
    \end{bmatrix}.
\end{align}
\endgroup
$\boldsymbol{D}$ has a special structure which allows us to efficiently compute $\boldsymbol{D}^{-1}$ via Schur complement methods. Note that we can re-write \cref{eq:ItI_D} as a block matrix with $2 \times 2$ blocks:
\begin{align}
\boldsymbol{D} &= \begin{bmatrix}
\boldsymbol{I}_{4 \nside \times 4 \nside} & \boldsymbol{D}_{12} \\
\boldsymbol{D}_{21} & \boldsymbol{I}_{4 \nside \times 4 \nside}
\end{bmatrix}.
\end{align}
This structure allows us to construct $\boldsymbol{W} = \boldsymbol{I}_{4 \nside \times 4 \nside}  - \boldsymbol{D}_{12} \boldsymbol{D}_{21}$, the Schur complement of the lower-right $\boldsymbol{I}_{4 \nside \times 4 \nside}$ block in $\boldsymbol{D}$; this is the only matrix we need to invert to compute $\boldsymbol{D}^{-1}$:
\begin{align}
    \boldsymbol{D}^{-1} &= \begin{bmatrix}
        \boldsymbol{W}^{-1} & -\boldsymbol{W}^{-1} \boldsymbol{D}_{12} \\
        - \boldsymbol{D}_{21}\boldsymbol{W}^{-1} & \boldsymbol{I}_{4 \nside \times 4 \nside} + \boldsymbol{D}_{21}\boldsymbol{W}^{-1}\boldsymbol{D}_{12}
    \end{bmatrix}. \label{eq:ItI_D_inv}
\end{align}
We use \cref{eq:ItI_D_inv} to compute $\boldsymbol{D}^{-1}$ and then construct the outputs of the merge stage using \cref{eq:merge_outputs_1,eq:merge_outputs_2}.

\section{Full algorithms for 3D problems using DtN matrices with a uniform discretization}
\label{appendix:3D_DtN_uniform}
In this section, we describe the details for the three-dimensional version of our method which merges DtN matrices. 
In this version of the algorithm, the outgoing boundary data $\boldsymbol{h}$ tabulates the outward-pointing normal derivative of the particular solution, and the incoming boundary data $\boldsymbol{g}$ tabulates the homogenous solution values restricted to patch boundaries. $\boldsymbol{T}$ is a Dirichlet-to-Neumann matrix.

As in the previous section, we use  $\boldsymbol{I}_{a\times a}$ to denote the identity matrix of shape $a \times a$ and $\boldsymbol{0}_{d}$ to denote a length-$d$ vector filled with $0$'s. 
We also use $\boldsymbol{0}$ to denote a matrix block filled with $0$'s, and assume the shape of this block can be inferred from its context.

\subsection{Local solve stage}
Recall from \cref{sec:discretization} that we use a tensor product of order-$p$ Chebyshev--Lobatto points to discretize the interior of each leaf. This results in a grid with $p^3$ discretization points, and $p^3 - (p-2)^3$ of these points lie on the boundary of the leaf. 
We use order-$q$ Gauss--Legendre points to discretize each side of the leaf's boundary, so there are $6q^2$ boundary points in total. 
Thus, to translate the information between the interior and boundary of each leaf, we need to compute spectral differentiation matrices and matrices interpolating between the $p^2$ Chebyshev and $q^2$ Gauss points on each face of the leaf. 
In particular, we need to precompute the following matrices:
\begin{itemize}
    \item $\boldsymbol{P}$, with shape $p^3 - (p-2)^3 \times 6q^2$, is the operator mapping data sampled on the Gauss boundary points to data sampled on the Chebyshev points located on the boundary of the leaf.
    This matrix is constructed using a barycentric Lagrange interpolation matrix mapping from Gauss to Chebyshev points on one face of the leaf; this interpolation matrix is repeated for the other sides. Rows corresponding to the Chebyshev points on the corners (edges) of the leaf average the contribution from the three (two) adjoining panels. 
    \item $\boldsymbol{Q}$ with shape  $6 q^2 \times p^3$, performs spectral differentiation on the $p^3$ Chebyshev points followed by interpolation to the Gauss boundary points.
    This matrix is formed by stacking the relevant rows of Chebyshev spectral differentiation matrices to form an operator which evaluates normal derivatives on the boundary Chebyshev points, and then composing this differentiation operator with a matrix formed from barycentric Lagrange interpolation matrix blocks. These interpolation matrices each map from one Chebyshev face to one Gauss face.
\end{itemize}
To work with $\boldsymbol{L}^{(i)}$, the discretization of the differential operator on leaf $i$, it is useful to identify $I_i$ and $I_e$, the sets of discretization points corresponding to the $(p-2)^3$ interior and $p^3 - (p-2)^3$ exterior Chebyshev points, respectively. Once the precomputed operators and index sets are correctly specified, \cref{alg:2D_DtN_local} can be re-used for the three-dimensional case.

\subsection{Merge stage}
\label{appendix:3D_DtN_merge}
\begin{figure}
    \centering
    \begin{subfigure}[b]{0.4\textwidth}
        \centering
        \includegraphics[width=\linewidth, trim={3cm 2cm 3cm 2cm}, clip]{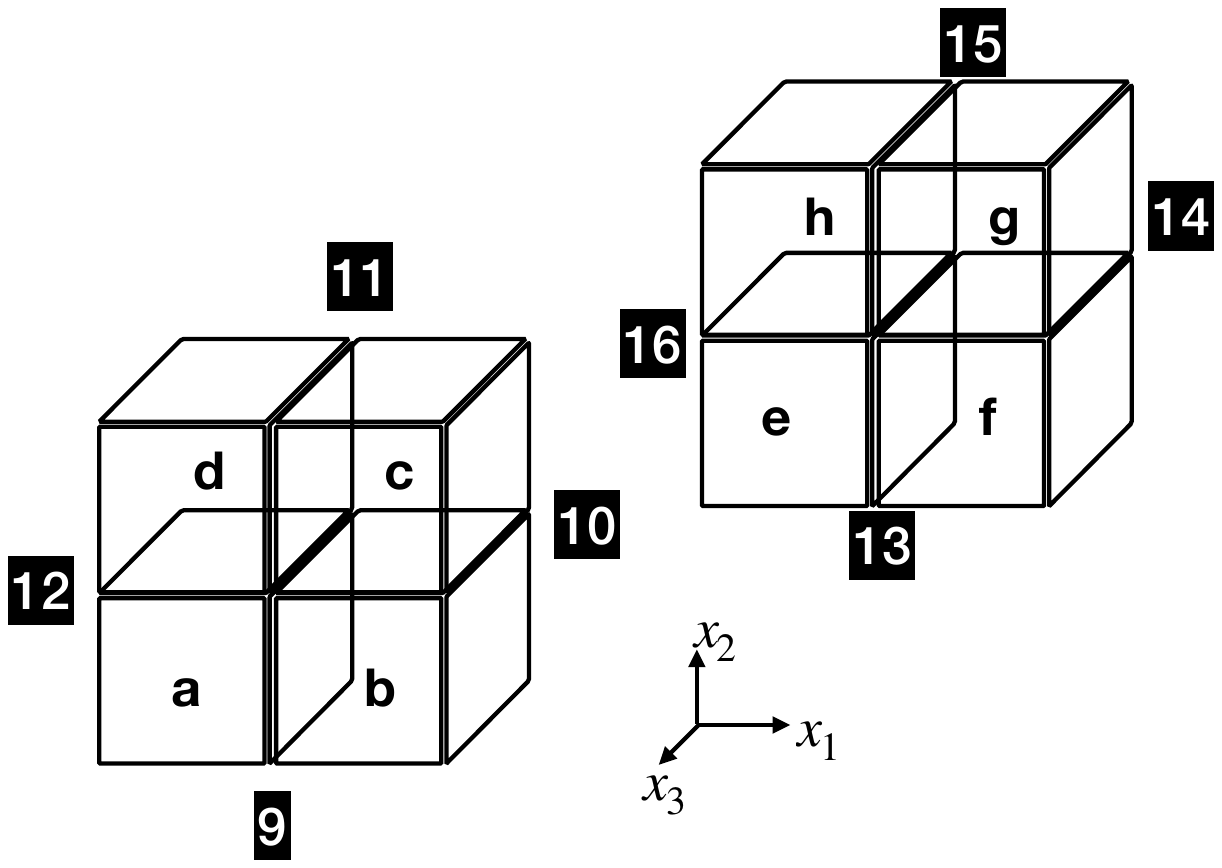}
    \end{subfigure}%
    \begin{subfigure}[b]{0.4\textwidth}
        \centering
        \includegraphics[width=\linewidth, trim={4cm 2cm 4cm 2cm}, clip]{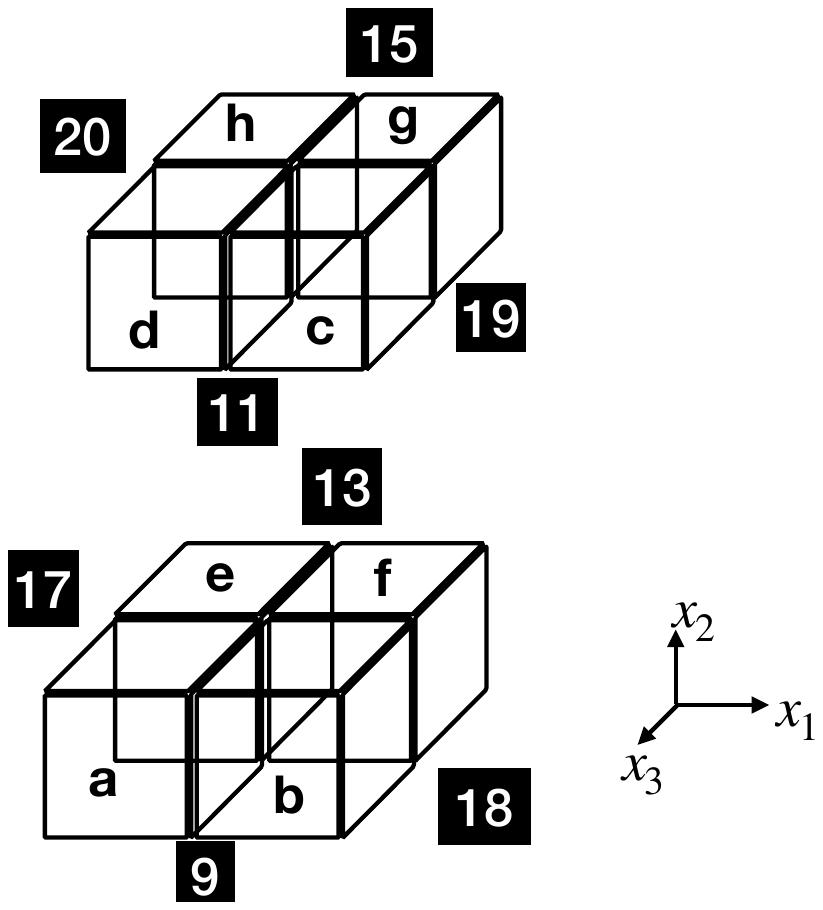}
    \end{subfigure}
\caption{Visualizing boundary elements $9$ through $20$ for three-dimensional merges.}
\label{fig:3D_merges}
\end{figure}
In the 3D merge stage, we are merging eight nodes $\Omega_a, \Omega_b, \Omega_c, \Omega_d, \Omega_e, \Omega_f, \Omega_g$, and $\Omega_h$, which have \emph{exterior} and \emph{interior} discretization points. We label the exterior boundary sections $1, 2, \hdots, 8$, and we label the interior boundary sections $9, 10, \hdots, 20$. See \cref{fig:3D_merges} for a diagram of the different boundary parts. 
Because the merge stage operates completely on data discretized using Gauss--Legendre panels, there are no discretization points at the corners or edges of nodes. This means each discretization point belongs to exactly one part of the boundary. 
During this stage of the algorithm, we will be indexing rows and columns of the Dirichlet-to-Neumann matrices according to boundary sections $1, \dots, 20$. 
For example, we use $\boldsymbol{T}^{(a)}_{1,9}$ to indicate the submatrix of node $a$'s DtN matrix which maps from boundary section $9$ to boundary section $1$. Suppose that each node has $\nside$ discretization points on each face. Then $\boldsymbol{T}^{(a)}_{1,9}$ will have shape $3 \nside \times \nside$.

Just as in the two-dimensional case, we use sets of constraints to solve for a mapping from given $\boldsymbol{g}_{\ext}$ to unknown $\boldsymbol{g}_{\intt}$, vectors tabulating the homogeneous solution along the boundary parts. First are constraints specifying that the solution to the PDE is continuous:
\begin{align}
    \boldsymbol{u}_{\ext} &= \boldsymbol{A} \boldsymbol{g}_{\ext} +  \boldsymbol{B}\boldsymbol{g}_{\intt} + \boldsymbol{h}^{(\child)}_{\ext}.
    \label{eq:3D_DtN_merge_1}
\end{align}
In this equation $\boldsymbol{u}_{\ext}$ is interpreted as the outward-pointing normal derivative of the solution to the PDE restricted to the merged nodes with boundary data specified by $\boldsymbol{g}_{\ext}$.
In this set of constraints, we define:
\begin{align}
    \boldsymbol{h}^{(\child)}_{\ext} &= \begin{bmatrix}
        \boldsymbol{h}^{(a)}_{1} \\
        \boldsymbol{h}^{(b)}_{2} \\
        \boldsymbol{h}^{(c)}_{3} \\
        \boldsymbol{h}^{(d)}_{4} \\
        \boldsymbol{h}^{(e)}_{5} \\
        \boldsymbol{h}^{(f)}_{6} \\
        \boldsymbol{h}^{(g)}_{7} \\
        \boldsymbol{h}^{(h)}_{8} \\
    \end{bmatrix}, \\
    \boldsymbol{A} &= \begin{bmatrix}
        \boldsymbol{T}^{(a)}_{1,1} & \bz & \bz & \bz & \bz & \bz & \bz & \bz \\ 
        \bz & \boldsymbol{T}^{(b)}_{2,2} & \bz & \bz & \bz & \bz & \bz & \bz \\ 
        \bz & \bz & \boldsymbol{T}^{(c)}_{3,3} & \bz & \bz & \bz & \bz & \bz \\
        \bz & \bz & \bz & \boldsymbol{T}^{(d)}_{4,4} & \bz & \bz & \bz & \bz \\
        \bz & \bz & \bz & \bz & \boldsymbol{T}^{(e)}_{5,5} & \bz & \bz & \bz \\
        \bz & \bz & \bz & \bz & \bz & \boldsymbol{T}^{(f)}_{6,6} & \bz & \bz \\
        \bz & \bz & \bz & \bz & \bz & \bz & \boldsymbol{T}^{(g)}_{7,7} & \bz \\
        \bz & \bz & \bz & \bz & \bz & \bz & \bz & \boldsymbol{T}^{(h)}_{8,8}
    \end{bmatrix},
\end{align}
\begin{equation}
    \resizebox{0.9\linewidth}{!}{
        $
        \boldsymbol{B} = \begin{bmatrix}
            \boldsymbol{T}^{(a)}_{1,9} & \bz & \bz & \boldsymbol{T}^{(a)}_{1,12} & \bz & \bz & \bz & \bz & \boldsymbol{T}^{(a)}_{1,17} & \bz & \bz & \bz \\
            \boldsymbol{T}^{(b)}_{2,9} & \boldsymbol{T}^{(b)}_{2,10} & \bz & \bz & \bz & \bz & \bz & \bz & \bz & \boldsymbol{T}^{(b)}_{2,18} & \bz & \bz  \\
            \bz & \boldsymbol{T}^{(c)}_{3,10} & \boldsymbol{T}^{(c)}_{3,11} & \bz & \bz & \bz & \bz & \bz & \bz & \bz & \boldsymbol{T}^{(c)}_{3,19} & \bz \\
            \bz & \bz & \boldsymbol{T}^{(d)}_{4,11} & \boldsymbol{T}^{(d)}_{4,12} & \bz & \bz & \bz & \bz & \bz & \bz & \bz & \boldsymbol{T}^{(d)}_{4,20} \\
            \bz & \bz & \bz & \bz & \boldsymbol{T}^{(e)}_{5,13} & \bz & \bz & \boldsymbol{T}^{(e)}_{5,16} & \boldsymbol{T}^{(e)}_{5,17} & \bz & \bz & \bz \\
            \bz & \bz & \bz & \bz & \boldsymbol{T}^{(f)}_{6,13} & \boldsymbol{T}^{(f)}_{6,14} & \bz & \bz & \bz & \boldsymbol{T}^{(f)}_{6,18} & \bz & \bz \\
            \bz & \bz & \bz & \bz & \bz  & \boldsymbol{T}^{(g)}_{7,14} & \boldsymbol{T}^{(g)}_{7,15} & \bz & \bz & \bz & \boldsymbol{T}^{(g)}_{7,19} & \bz \\
            \bz & \bz & \bz & \bz & \bz & \bz & \boldsymbol{T}^{(h)}_{8,15} & \boldsymbol{T}^{(h)}_{8,16} & \bz & \bz & \bz & \boldsymbol{T}^{(h)}_{8,20}
        \end{bmatrix}.
        $
    }
\end{equation}

As in the two-dimensional case, we enforce constraints that ensure the normal derivatives from neighboring nodes sum to zero, which gives us a system of constraints:
\begin{align}
    \bz_{12 \nside} &= \boldsymbol{C} \boldsymbol{g}_{\ext} +  \boldsymbol{D}\boldsymbol{g}_{\intt} + \boldsymbol{h}^{(\child)}_{\intt},
    \label{eq:3D_DtN_merge_2}
\end{align}
where we define
\begin{align}
    \boldsymbol{h}_{\intt}^{(\child)} &= \begin{bmatrix}
        \boldsymbol{h}^{(a)}_{9} + \boldsymbol{h}^{(b)}_{9} \\
        \boldsymbol{h}^{(b)}_{10} + \boldsymbol{h}^{(c)}_{10} \\
        \boldsymbol{h}^{(c)}_{11} + \boldsymbol{h}^{(d)}_{11} \\
        \boldsymbol{h}^{(d)}_{12} + \boldsymbol{h}^{(a)}_{12} \\
        \boldsymbol{h}^{(e)}_{13} + \boldsymbol{h}^{(f)}_{13} \\
        \boldsymbol{h}^{(f)}_{14} + \boldsymbol{h}^{(g)}_{14} \\
        \boldsymbol{h}^{(g)}_{15} + \boldsymbol{h}^{(h)}_{15} \\
        \boldsymbol{h}^{(h)}_{16} + \boldsymbol{h}^{(e)}_{16} \\
        \boldsymbol{h}^{(a)}_{17} + \boldsymbol{h}^{(e)}_{17} \\
        \boldsymbol{h}^{(b)}_{18} + \boldsymbol{h}^{(f)}_{18} \\
        \boldsymbol{h}^{(c)}_{19} + \boldsymbol{h}^{(g)}_{19} \\
        \boldsymbol{h}^{(d)}_{20} + \boldsymbol{h}^{(h)}_{20}
    \end{bmatrix}, \\
    \boldsymbol{C} &= \begin{bmatrix}
        \boldsymbol{T}^{(a)}_{9,1} & \boldsymbol{T}^{(b)}_{9,2} & \bz & \bz & \bz & \bz & \bz & \bz \\
        \bz & \boldsymbol{T}^{(b)}_{10,2} & \boldsymbol{T}^{(c)}_{10,3} & \bz & \bz & \bz & \bz & \bz \\
        \bz & \bz & \boldsymbol{T}^{(c)}_{11,3} & \boldsymbol{T}^{(d)}_{11,4} & \bz & \bz & \bz & \bz \\
        \boldsymbol{T}^{(a)}_{12,1} & \bz & \bz & \boldsymbol{T}^{(d)}_{12,4} & \bz & \bz & \bz & \bz \\
        \bz & \bz & \bz & \bz & \boldsymbol{T}^{(e)}_{13,5} & \boldsymbol{T}^{(f)}_{13,6} & \bz & \bz \\
        \bz & \bz & \bz & \bz & \bz & \boldsymbol{T}^{(f)}_{14,6} & \boldsymbol{T}^{(g)}_{14,7} & \bz \\
        \bz & \bz & \bz & \bz & \bz & \bz & \boldsymbol{T}^{(g)}_{15,7} & \boldsymbol{T}^{(h)}_{15,8} \\
        \bz & \bz & \bz & \bz & \boldsymbol{T}^{(e)}_{16,5} & \bz & \bz & \boldsymbol{T}^{(h)}_{16,8} \\
        \boldsymbol{T}^{(a)}_{17,1} & \bz & \bz & \bz & \boldsymbol{T}^{(e)}_{17,5} & \bz & \bz & \bz \\
        \bz & \boldsymbol{T}^{(b)}_{18,2} & \bz & \bz & \bz & \boldsymbol{T}^{(f)}_{18,6} & \bz & \bz \\
        \bz & \bz & \boldsymbol{T}^{(c)}_{19,3} & \bz & \bz & \bz & \boldsymbol{T}^{(g)}_{19,7} & \bz \\
        \bz & \bz & \bz & \boldsymbol{T}^{(d)}_{20,4} & \bz & \bz & \bz & \boldsymbol{T}^{(h)}_{20,8}
    \end{bmatrix},
\end{align}
\begin{equation}
    \resizebox{0.91\linewidth}{!}{
        $
    \boldsymbol{D} = \begin{bmatrix}
        \boldsymbol{T}^{(a)}_{9,9} + \boldsymbol{T}^{(b)}_{9,9} & \boldsymbol{T}^{(b)}_{9, 10} & \bz & \boldsymbol{T}^{(a)}_{9, 12} & \bz & \bz & \bz & \bz & \boldsymbol{T}^{(a)}_{9, 17} & \boldsymbol{T}^{(b)}_{9, 18} & \bz & \bz \\
        \boldsymbol{T}^{(b)}_{10,9} & \boldsymbol{T}^{(b)}_{10,10} + \boldsymbol{T}^{(c)}_{10,10} & \boldsymbol{T}^{(c)}_{10,11} & \bz & \bz & \bz & \bz & \bz & \bz & \boldsymbol{T}^{(b)}_{10, 18} & \boldsymbol{T}^{(c)}_{10,19} & \bz \\
        \bz &  \boldsymbol{T}^{(c)}_{11,10} & \boldsymbol{T}^{(c)}_{11,11} + \boldsymbol{T}^{(d)}_{11,11} & \boldsymbol{T}^{(d)}_{11,12} & \bz & \bz & \bz & \bz & \bz & \bz &  \boldsymbol{T}^{(c)}_{11,19} & \boldsymbol{T}^{(d)}_{11,20} \\
        \boldsymbol{T}^{(a)}_{12,9} & \bz & \boldsymbol{T}^{(d)}_{12,11} & \boldsymbol{T}^{(d)}_{12,12} + \boldsymbol{T}^{(a)}_{12,12} & \bz & \bz & \bz & \bz & \boldsymbol{T}^{(a)}_{12,17} & \bz & \bz & \boldsymbol{T}^{(d)}_{12,20} \\
        \bz & \bz & \bz & \bz & \boldsymbol{T}^{(e)}_{13,13} + \boldsymbol{T}^{(f)}_{13,13} & \boldsymbol{T}^{(f)}_{13,14} & \bz & \boldsymbol{T}^{(e)}_{13,16} & \boldsymbol{T}^{(e)}_{13,17} & \boldsymbol{T}^{(f)}_{13,18} & \bz & \bz \\
        \bz & \bz & \bz & \bz & \boldsymbol{T}^{(f)}_{14,13} & \boldsymbol{T}^{(f)}_{14,14} + \boldsymbol{T}^{(g)}_{14,14} & \boldsymbol{T}^{(g)}_{14,15} & \bz & \bz & \boldsymbol{T}^{(f)}_{14,18} & \boldsymbol{T}^{(g)}_{14,19} & \bz \\
        \bz & \bz & \bz & \bz & \bz & \boldsymbol{T}^{(g)}_{15,14} & \boldsymbol{T}^{(g)}_{15,15} + \boldsymbol{T}^{(h)}_{15,15} & \boldsymbol{T}^{(h)}_{15,16} & \bz & \bz & \boldsymbol{T}^{(g)}_{15,19} & \boldsymbol{T}^{(h)}_{15,20} \\
        \bz & \bz & \bz & \bz & \boldsymbol{T}^{(e)}_{16,13} & \bz & \boldsymbol{T}^{(h)}_{16,15} & \boldsymbol{T}^{(h)}_{16,16} + \boldsymbol{T}^{(e)}_{16,16} & \boldsymbol{T}^{(e)}_{16,17} & \bz & \bz & \boldsymbol{T}^{(h)}_{16,20} \\
        \boldsymbol{T}^{(a)}_{17,9} & \bz & \bz & \boldsymbol{T}^{(a)}_{17,12} & \boldsymbol{T}^{(e)}_{17,13} & \bz & \bz & \boldsymbol{T}^{(e)}_{17,16} & \boldsymbol{T}^{(a)}_{17,17} + \boldsymbol{T}^{(e)}_{17,17} & \bz & \bz & \bz \\
        \boldsymbol{T}^{(b)}_{18,9} & \boldsymbol{T}^{(b)}_{18,10} & \bz & \bz & \boldsymbol{T}^{(f)}_{18,13} & \boldsymbol{T}^{(f)}_{18,14} & \bz & \bz & \bz & \boldsymbol{T}^{(b)}_{18,18} + \boldsymbol{T}^{(f)}_{18,18} & \bz & \bz \\
        \bz & \boldsymbol{T}^{(c)}_{19,10} & \boldsymbol{T}^{(c)}_{19,11} & \bz & \bz & \boldsymbol{T}^{(g)}_{19,14} & \boldsymbol{T}^{(g)}_{19,15} & \bz & \bz & \bz & \boldsymbol{T}^{(c)}_{19,19} + \boldsymbol{T}^{(g)}_{19,19} & \bz \\
        \bz & \bz  & \boldsymbol{T}^{(d)}_{20,11} & \boldsymbol{T}^{(d)}_{20,12} & \bz &  \bz & \boldsymbol{T}^{(h)}_{20,15} & \boldsymbol{T}^{(h)}_{20,16} & \bz & \bz & \bz & \boldsymbol{T}^{(d)}_{20,20} + \boldsymbol{T}^{(h)}_{20,20} \\
    \end{bmatrix}.
    $
    }
\end{equation}
Now that the matrices and vectors are defined, we can construct the linear system in \cref{eq:merge_lin_system} and compute the merged data.

\section{Details for the inverse scattering experiment}
\label{appendix:inverse_scattering}
\revone{We consider a forward model which composes $\cF$, which maps a scattering potential to evaluations of a scattered wave, and $\cB$, which maps coefficients of a sine series to scattering potentials:}
\begin{align}
    \revone{\left( \cF \circ \cB \right) (\theta)_j} &\revone{= u_\theta(x^{(j)}); }& \revone{u_\theta \text{ solves Equation }  \mbox{\ref{eq:forward_scattering_problem}} \text{ with }q=q_\theta;} \\
    \revone{q_\theta(x)} & \revone{= \cB(\theta) = \sum_{b_j \in B_\gamma }\theta_j b_j(x);}
\end{align}
\revone{where the basis $B_\gamma$ is specified by \mbox{\cref{eq:basis_def}}. In this experiment, we use a domain $\Omega = [-1, 1]$, and we use $100$ evaluation points equispaced on a ring of radius $5$:}
\begin{align}
    \label{eq:disc_pts}
    \revone{x^{(j)} = \left( 5 \sin\left( \frac{2 \pi j}{100} \right), 5 \cos \left( \frac{2 \pi j}{100} \right) \right), \qquad j=0,\dots, 99.}
\end{align}
\revone{To solve an inverse problem, we are interested in computing the Fr\'echet derivative of $\cF \circ \cB$ centered at $\theta$; we call this object $J[\theta]$. By the chain rule, we can decompose this Fr\'echet derivative} 
\begin{align*}
    \revone{J[\theta] = J_{\cF}[q_\theta] \circ J_\cB[\theta].}
\end{align*}
\revone{The basis transformation $\cB$ is linear, so the action of $J_\cB[\theta]$ can be computed with standard sine and adjoint sine transforms. The following section will describe how we compute the action of $J_\cF [q_\theta]$.}
\subsection{Defining the action of the Fr\'echet derivative}
\label{appendix:actions}
\revone{
\mbox{\citet{borges_high_2016}} describe the Fr\'echet derivative $J_\cF [q_\theta]$ the action $J_\cF [q_\theta] v$ and $J_\cF [q_\theta]^* f$. The action of the derivative and its adjoint can be described by the solution of elliptic partial differential equations. We re-state these results in this section.}
\begin{theorem}[\revone{Theorem 3.1 of \mbox{\citet{borges_high_2016}}}]
    \label{thm:derivative}
    \revone{Let $u_\theta$ solve \mbox{\cref{eq:forward_scattering_problem}} with scattering potential $q=q_\theta$. Let $w$ solve }
    \begin{equation}
\left\{
    \begin{aligned}
        \revone{\Delta w(x) + k^2(1 + q(x)) w(x)} &\revone{= k^2 v(x) \left( u_\theta(x) + e^{ik\langle \hat{s}, x \rangle} \right),} && \quad \revone{x \in [-1, 1]^2,} \\
        \revone{\sqrt{r} \left( \tfrac{\partial w}{\partial r} - ikw \right)} &\revone{\to 0,} && \revone{\quad r = \|x \|_2 \to \infty.}
    \end{aligned}
\right.
    \label{eq:forward_derivative}
\end{equation}
\revone{Then}
\begin{align*}
    \revone{\left( J_\cF [q_\theta] v \right)_j} &\revone{= w(x^{(j)}).}
\end{align*}
\end{theorem}

\begin{theorem}[\revone{Theorem 3.2 of \mbox{\citet{borges_high_2016}}}]
    \label{thm:adjoint}
    \revone{
    Let $u_\theta$ solve \mbox{\cref{eq:forward_scattering_problem}} with scattering potential $q=q_\theta$. Let $f$ denote a singular charge distribution supported on the evaluation points $\{ x^{(j)} \}_{j=0, \dots, 99}$ viewed as a generalized function in $\R^2$. Let $w$ solve }
    \begin{equation}        
    \left\{
        \begin{aligned}
            \revone{\Delta w(x) + k^2(1 + q(x)) w(x)} &\revone{= k^2 f(x),} && \quad \revone{x \in \R^2,} \\
            \revone{\sqrt{r} \left( \tfrac{\partial w}{\partial r} + ikw \right)} &\revone{\to 0,} && \revone{\quad r = \|x \|_2 \to \infty.}
        \end{aligned}
    \right.
    \end{equation}
    \revone{Then} 
    \begin{align*}
    \revone{\left( J_\cF [q_\theta]^* f \right)(x)} &\revone{= w(x)\overline{u_\theta(x) + e^{ik\langle \hat{s}, x \rangle}}.}
\end{align*}
\end{theorem}

\subsection{Experiment details}
\label{app:adjoints_experimental_details}
\revone{
In this section, we describe the experimental setting used to generate \mbox{\cref{fig:autodiff_accuracy}}. In these experiments, we use $k=20$ and basis $B_\gamma$ with $\gamma = 25$, which gives us $N_\theta = 465$ basis coefficients. 
We compute a coefficient vector $\theta$ by projecting the scattering potential \mbox{\cref{eq:q_sum_gauss}} onto $B_\gamma$.
}

\revone{To measure the accuracy of the outputs of JAX Jacobian-vector products, we generate a random coefficient vector $\delta$ distributed i.i.d. Gaussian and compute $J[\theta] \delta$. We compute this Jacobian-vector product for a range of discretization sizes, holding $p=16$ constant and varying $L=1,\hdots, 5$.
We compare this with the action of the Fr\'echet derivative, which is computed by first computing $J_\cB[\theta] \delta$ and then applying $J_\cF[q_\theta]$ from \revone{\cref{thm:derivative}}, which is computed using parameters $p=16$ and $L=5$. We measure the relative $\ell_\infty$ error between the outputs of JAX Jacobian-vector products and the the action of the Fr\'echet derivative.}

\revone{To measure the accuracy of the outputs of JAX vector-Jacobian products, we generate a random perturbation $f \in \C^{100}$. The real and imaginary parts of each component of $f$ are distributed i.i.d. Gaussian. We compute this vector-Jacobian product for a range of discretization sizes, again holding $p=16$ constant and varying $L=1, \hdots, 5$. We compare this with the action of the Fr\'echet derivative which is computed by first evaluating $J_\cF[q_\theta]^*f$ via \mbox{\cref{thm:adjoint}}, which is computed using parameters $p=16$ and $L=5$, and then applying $J_\cB[\theta]^*$. We measure the relative $\ell_\infty$ error betwen the outputs of JAX vector-Jacobian products and the action of the Fr\'echet derivative.}

\revone{We note the choice of evaluation points exterior to the computational domain (\mbox{\cref{eq:disc_pts}}) is not necessary for the convergence of automatic differentiation; we choose these evaluation points to ensure numerical stability when computing $J_{\cF}[q_\theta]^* f$. In preliminary experiments, we observed accurate vector-Jacobian products when the evaluation points were located at the HPS discretization points. In this setting, the inputs of the vector-Jacobian product routine must be scaled by the appropriate quadrature weights.
}

\section{Additional timing results}
\label{sec:addl_results}
\revtwo{In this appendix, we extend the results shown in \mbox{\cref{fig:recomputation_results}} to include our subtree recomputation method evaluated with different subtree depths. 
As a hueristic, we choose to use the subtree recomputation depth to be the maximum depth where all computations for \mbox{\cref{alg:local_solve,alg:merge}} can fit on a single GPU. In \mbox{\cref{fig:recomputation_results}}, we consider the DtN version of the method with polynomial order $p=16$; this results in subtree depth $7$. In \mbox{\cref{tab:extra}}, we show the results of choosing different subtree depth parameters. We measure the runtime for large problem sizes $L=8$ and $9$; recomputation is necessary for these problem sizes.}

\revtwo{As we decrease the subtree depth, we see runtimes increase, as more transfers between the GPU and host RAM are required.}

\begin{table}
    \centering
\begin{tabular}{lrrrr}
    \revtwo{\textbf{Subtree depth}} & \revtwo{$L$} & \revtwo{$N$} & \revtwo{\textbf{Runtime (s)}} & \revtwo{\textbf{\% of Peak FLOPS}} \\
    \midrule
    \revtwo{5} & \revtwo{8} & \revtwo{$16{,}777{,}216$} & \revtwo{$5.86$} & \revtwo{$6.68\%$}  \\
    \revtwo{5} & \revtwo{9} & \revtwo{$67{,}108{,}864$} & \revtwo{$24.50$} & \revtwo{$10.86\%$}\\
    \midrule
    \revtwo{6} & \revtwo{8} & \revtwo{$16{,}777{,}216$} & \revtwo{$4.35$} & \revtwo{$10.58\%$}  \\
    \revtwo{6} & \revtwo{9} & \revtwo{$67{,}108{,}864$} & \revtwo{$18.85$} & \revtwo{$15.59\%$}\\
    \midrule
    \revtwo{7} & \revtwo{8} & \revtwo{$16{,}777{,}216$} & \revtwo{$4.02$} & \revtwo{$14.86\%$}  \\
    \revtwo{7} & \revtwo{9} & \revtwo{$67{,}108{,}864$} & \revtwo{$17.43$} & \revtwo{$20.01\%$}\\
    \bottomrule
\end{tabular}
\caption{\revtwo{Evaluating the effect of the subtree height parameter.}}
\label{tab:extra}
\end{table}

\end{document}